\documentclass{wstmp}

\usepackage{bm}
\usepackage{enumerate}

\newcommand{\RR}{\mathbb{R}}
\newcommand{\ZZ}{\mathbb{Z}}
\newcommand{\NN}{\mathbb{N}}
\newcommand{\QQ}{\mathbb{Q}}

\newcommand{\cB}{{\mathcal B}}
\newcommand{\cC}{{\mathcal C}}
\newcommand{\cE}{{\mathcal E}}
\newcommand{\cK}{{\mathcal K}}
\newcommand{\cL}{{\mathcal L}}
\newcommand{\cM}{{\mathcal M}}
\newcommand{\cP}{{\mathcal P}}
\newcommand{\cS}{{\mathcal S}}
\newcommand{\cT}{{\mathcal T}}

\newcommand{\tri}[3]{\langle #1, #2 \, | \, #3 \rangle}

\newcommand{\convcone}{\RR^+}

\newcommand{\bra}[1]{\{ #1 \}}
\newcommand{\indtri}[4]{#1 \indep #2 \, | \, #3 \, [#4] }
\newcommand{\indtr}[3]{#1 \indep #2 \, | \, #3}
\newcommand{\indep}{\mathop{\perp\!\!\!\perp}}

\newcommand{\cfg}{\mathcal{U}}

\newcommand{\precw}{<}
\newcommand{\succw}{>}

\newcommand{\cG}{{\mathcal G}}
\newcommand{\utri}[3]{u_{\langle #1, #2 \, | \, #3 \rangle}}
\newcommand{\abs}[1]{\left| #1 \right|}
\newcommand{\abC}[1]{a_{#1}b_{#1}C_{#1}}
\newcommand{\abCd}[1]{a_{#1}'b_{#1}'C_{#1}'}

\begin{document}

\title{CONES OF ELEMENTARY IMSETS AND SUPERMODULAR FUNCTIONS: A REVIEW AND SOME NEW RESULTS}
\author{Takuya Kashimura}
\address{Department of 
Mathematical Informatics,\\
Graduate School of Information Science and Technology, \\
University of Tokyo}
\author{Tomonari Sei}
\address{Department of Mathematics, Keio University}
\author{Akimichi Takemura}
\address{Department of 
Mathematical Informatics,\\
Graduate School of Information Science and Technology, \\
University of Tokyo\\
JST, CREST}
\author{Kentaro Tanaka}
\address{Department of Industrial Engineering and Management,\\
Graduate School of Decision Science and Technology,\\
Tokyo Institute of Technology}



\begin{abstract}
In this paper we give a review of the method of imsets introduced by 
Studen\'y\cite{stu2005} from a geometric point of view.  
Elementary imsets span a polyhedral cone and its dual cone is the cone
of supermodular functions.  We review basic facts on
the structure of these cones.  
Then we derive some new results on the following topics: 
i) extreme rays of the cone of standardized supermodular functions,
ii) faces of the cones, 
iii) small relations among elementary imsets,
and 
iv) some computational results on Markov basis for the toric ideal defined by elementary imsets.
\end{abstract}

\section{Introduction}
\label{sec:introduction}

The method of imsets by Studen\'y\cite{stu2005} provides a very
powerful algebraic method for describing conditional independence
relations under a probability measure.  Rules for deriving 
conditional independence relations are translated into 
relations among integer vectors called imsets.  
Hence
many properties of conditional independence relations
can be conveniently interpreted from a geometric viewpoint.
In recent papers Studen\'y and his collaborators
\cite{studeny-vomlel-hemmecke-2010,studeny-vomlel-2011}
further develop geometric methods for learning Bayesian networks.
In this paper we are more concerned on basic geometric
properties of imsets, in particular from the viewpoint of lattice  bases
and Markov bases for the configuration of elementary imsets.

The cone of supermodular functions, which we call the {\it supermodular cone},
is defined by a set of
linear inequalities and effective  
inequalities
correspond to elementary imsets.  Hence the $H$-representation
(cf.\ Gr\"unbaum\cite{grunbaum-2003}) 
of the supermodular cone is explicitly given.  The cone generated
by the elementary imsets, which we call the {\it imset cone}, 
is the dual to the supermodular cone and its 
set of extreme rays is given by elementary imsets.
Hence, in the dual sense, the $V$-representation of the imset cone
is given.  

From an algorithmic viewpoint of convex geometry, 
the $V$-representation of the supermodular cone, or equivalently the
$H$-representation of the imset cone, are hard to compute and
characterize.  
Therefore 
general results on the facets of
the imset cone, or equivalently 
the ``extreme rays'' of the supermodular cone
are important.
Since the supermodular cone contains a linear subspace consisting
of modular functions, we consider extreme rays of the cone of standardized
supermodular functions (see \eqref{eq:super-modular-standardization} below
for standardization).
These extreme rays are called skeletal supermodular functions in
Studen\'y\cite{stu2005}. 

Although the complete description of the supermodular cone
and the imset cone is very difficult, some faces of these
cones can be studied in detail.  In particular a face
corresponding to a semi-elementary imset seems to have a simpler
structure than other faces.  

The organization of the paper is as follows.  In Sections
\ref{sec:basic-facts} and \ref{sec:multiinformation} we present a
review of the method of imsets from a geometric point of
view.  In Section \ref{sec:basic-facts} we summarize basic facts on
the cones of supermodular functions and semi-elementary imsets.  In Section
\ref{sec:multiinformation} we review derivations of conditional
independence statements based on multiinformation and imsets.
In Section \ref{sec:extreme-rays} we give some results
on extreme rays of the cone of standardized supermodular functions.
In Section \ref{sec:semi-elementary-face} we discuss
properties of a face corresponding to a semi-elementary imset.
In Section \ref{sec:small-linear-relations} we characterize
relations among a small number of elementary imsets.
Finally in Section \ref{sec:markov-basis} we discuss some
computational results on Markov basis for the toric ideal 
defined by the configuration of elementary imsets.

\section{Basic facts on supermodular functions and imsets}
\label{sec:basic-facts}

As in {Kashimura and Takemura\cite{kashimura-takemura-2011}}, 
we use the notation and definitions from Studen\'y\cite{stu2005}.
Let $N$ be a finite set  and let $\mathcal{P}(N) = \bra{A: A \subseteq N}$ 
denote its power set.
In this paper $A\subseteq B$ means that $A$ is a subset of $B$
and $A\subset B$ means that $A$ is a proper subset of $B$.
$|A|$ denotes the cardinality of $A$.
For notational convenience, 
we write the union $A \cup B$ as $AB$.
A singleton set $\{i\}$ is simply written as $i$.
$\RR$, $\RR^+$,
$\QQ$, $\QQ^+$,
$\ZZ$, $\ZZ^+$,
$\NN$, 
denote the sets of reals, non-negative reals, rationals, non-negative rationals,
integers, non-negative integers and positive integers, 
respectively.

$f :\cP(N)\rightarrow \RR$ is called  {\em supermodular} if
\begin{equation}
\label{eq:super-modular-function}
f(EF)+ f(E\cap F)\ge f(E)+f(F), \qquad \forall E,F\subseteq N.
\end{equation}
The set of supermodular functions over $N$ is denoted by $\cK(N)$.
$f$ is submodular if $-f$ is supermodular.

$f$ is  {\em modular} if
it is both supermodular and submodular, i.e.\ 
\begin{equation}
\label{eq:modular-function}
f(EF)+ f(E\cap F)=f(E)+f(F), \qquad \forall E,F\subseteq N.
\end{equation}
$\cL(N)=\cK(N)\cap (-\cK(N))$ denotes the set of modular functions over $N$.
A modular function $f$ is like a discrete (signed) measure.
Indeed if $f(\emptyset)=0$, then by taking disjoint $E$ and $F$ 
in \eqref{eq:modular-function}, we see
that $f$ is a measure and hence $f$ can be written as  $f(E)=\sum_{e\in E} f(e)$.
Without the restriction of $f(\emptyset)=0$, we have
$
f(E)-f(\emptyset)=\sum_{e\in E} (f(e)-f(\emptyset)), 
$
which can also be written as
\begin{equation}
\label{eq:modular}
f(E)=\lambda_\emptyset + \sum_{e\in E} \lambda_e, \qquad \lambda_\emptyset = f(\emptyset), \ \lambda_e = f(e)-f(\emptyset), e\in N.
\end{equation}
This shows that the dimension of the linear space $\cL(N)$ is $|N|+1$
and a basis of $\cL(N)$ is given by the following $|N|+1$ functions:
\begin{equation}
\label{eq:modular-basis}
f_\emptyset(E)\equiv 1, \quad f_e(E)=1_{e\in\bullet}(E)=1_{e\in E}, \quad e\in N, 
\end{equation}
where $1_{e\in\bullet}(E)
=1_{\{e\}\subseteq\bullet}(E)=1_{e\in E}$ is the indicator function
\begin{equation}
\label{eq:indicator}
1_{e\in\bullet} (E)= \begin{cases} 1 & e\in E, \\
                       0 & \text{otherwise}.
\end{cases}
\end{equation}

Given a supermodular function $f$, define 
a modular function $f_L$ and a supermodular function $\bar f$ by 
\begin{align}
f_L(E)&= \lambda_\emptyset + \sum_{e\in E} \lambda_e, \quad \lambda_\emptyset = f(\emptyset),
\quad \lambda_e = f(e)-f(\emptyset),  \nonumber \\
\bar f(E)&=f(E)-f_L(E).
\label{eq:super-modular-standardization}
\end{align}
Then  $\bar f$ is supermodular and $\bar f(E)=0$ for $|E|\le 1$.
$\bar f$ is often a preferred standardization of $f$.
By induction on the cardinality of $|E|$, it is easy to show that
any standardized supermodular function $f$ is non-negative and
non-decreasing:
\[
f(E)\ge 0, \quad E\subseteq F \Rightarrow f(E)\le f(F), \quad E,F\subseteq N.
\]

Let $N=\{1,\dots,|N|\}$.  We can identify $f:\cP(N)\rightarrow \RR$ with a vector in
$\RR^{| \cP(N)|}=\RR^{2^{|N|}}$ by listing its values:
\[
\big(f(\emptyset), f(1),\dots,f(|N|), f({\{1,2\}}), \dots, f(N)\big).
\]
With this identification the set of supermodular functions $\cK(N)$ is
a polyhedral cone in $\RR^{2^{|N|}}$, because it is defined
by the set of linear inequalities in \eqref{eq:super-modular-function}.
By the definition of $\cL(N)=\cK(N)\cap (-\cK(N))$, $\cL(N)$ 
is the largest linear subspace contained in $\cK(N)$.

For $A\subseteq N$,  $\delta_A: \cP(N)\rightarrow \RR$
is defined as
\[
\delta_A(S)=\begin{cases} 1 & \text{if } S=A, \\
   0 & \text{otherwise}.
 \end{cases}
\]

For pairwise disjoint subsets, $A,B,C \subseteq N$,
we write this triplet by $\tri{A}{B}{C}$, and 
the set of all disjoint triplets $\tri{A}{B}{C}$ over $N$ by $\mathcal{T}(N)$.
Unless otherwise stated, we assume that $A,B$ are non-empty. 
On the other hand $C$ may well be an empty set.
For a triplet $\tri{A}{B}{C} \in \mathcal{T}(N)$, a {\it semi-elementary imset} $u_{\tri{A}{B}{C}} : {\cal P}(N) \rightarrow \RR$ 
is defined as 
\begin{align}
        u_{\tri{A}{B}{C}} = \delta_{ABC} + \delta_C - \delta_{AC} - \delta_{BC}. 
\label{eq:semi-elementary}
\end{align}
If $A=\emptyset$ or $B=\emptyset$, then $u_{\tri{A}{B}{C}}$ is the zero imset.
We consider $u_{\tri{A}{B}{C}}$ as a $2^{|N|}$-dimensional integer vector
with two elements equal to $1$ (at $ABC$ and $C$) and two elements equal to 
$-1$ (at $AC$ and $BC$).
If $A = a$ and $B = b$ are singletons, 
the imset $u_{\tri{a}{b}{C}}$ is called {\it elementary}.
The set of all elementary imsets is denoted by $\cE(N)$.  Note that
the number of elementary imsets is given by 
\[
|\cE(N)| = \binom{|N|}{2} 2^{|N|-2}=|N|(|N|-1)2^{|N|-3}.
\]


\newcommand{\bs}{\begin{sideways}$}
\newcommand{\es}{$\end{sideways}}
\begin{table}[thbp]
\tbl{Configuration of elementary imsets $\cfg_N$, $|N|=4$.}{\hspace*{8cm}}
\label{tab:configuration}
{\scriptsize
\setlength{\tabcolsep}{3pt}
\begin{tabular}{c|cccccccccccccccccccccccc}
&
\bs \tri{a}{b}{cd} \es &
\bs \tri{a}{c}{bd} \es &
\bs \tri{a}{d}{bc} \es &
\bs \tri{b}{c}{ad} \es &
\bs \tri{b}{d}{ac} \es &
\bs \tri{c}{d}{ab} \es &
\bs \tri{b}{c}{d} \es &
\bs \tri{a}{c}{d} \es &
\bs \tri{a}{b}{d} \es &
\bs \tri{b}{d}{c} \es &
\bs \tri{a}{d}{c} \es &
\bs \tri{a}{b}{c} \es &
\bs \tri{c}{d}{b} \es &
\bs \tri{a}{d}{b} \es &
\bs \tri{a}{c}{b} \es &
\bs \tri{c}{d}{a} \es &
\bs \tri{b}{d}{a} \es &
\bs \tri{b}{c}{a} \es &
\bs \tri{c}{d}{\emptyset} \es &
\bs \tri{b}{d}{\emptyset} \es &
\bs \tri{a}{d}{\emptyset} \es &
\bs \tri{b}{c}{\emptyset} \es &
\bs \tri{a}{c}{\emptyset} \es &
\bs \tri{a}{b}{\emptyset} \es 
\\ \hline
 $N$         &  1 &  1 &  1 &  1 &  1 &  1 &  0 &  0 &  0 &  0 &  0 &  0 &  0 &  0 &  0 &  0 &  0 &  0 &  0 &  0 &  0 &  0 &  0 &  0 \\
 $bcd$       & -1 & -1 & -1 &  0 &  0 &  0 &  1 &  0 &  0 &  1 &  0 &  0 &  1 &  0 &  0 &  0 &  0 &  0 &  0 &  0 &  0 &  0 &  0 &  0 \\
 $acd$       & -1 &  0 &  0 & -1 & -1 &  0 &  0 &  1 &  0 &  0 &  1 &  0 &  0 &  0 &  0 &  1 &  0 &  0 &  0 &  0 &  0 &  0 &  0 &  0 \\
 $abd$       &  0 & -1 &  0 & -1 &  0 & -1 &  0 &  0 &  1 &  0 &  0 &  0 &  0 &  1 &  0 &  0 &  1 &  0 &  0 &  0 &  0 &  0 &  0 &  0 \\
 $abc$       &  0 &  0 & -1 &  0 & -1 & -1 &  0 &  0 &  0 &  0 &  0 &  1 &  0 &  0 &  1 &  0 &  0 &  1 &  0 &  0 &  0 &  0 &  0 &  0 \\
 $cd$        &  1 &  0 &  0 &  0 &  0 &  0 & -1 & -1 &  0 & -1 & -1 &  0 &  0 &  0 &  0 &  0 &  0 &  0 &  1 &  0 &  0 &  0 &  0 &  0 \\
 $bd$        &  0 &  1 &  0 &  0 &  0 &  0 & -1 &  0 & -1 &  0 &  0 &  0 & -1 & -1 &  0 &  0 &  0 &  0 &  0 &  1 &  0 &  0 &  0 &  0 \\
 $bc$        &  0 &  0 &  1 &  0 &  0 &  0 &  0 &  0 &  0 & -1 &  0 & -1 & -1 &  0 & -1 &  0 &  0 &  0 &  0 &  0 &  0 &  1 &  0 &  0 \\
 $ad$        &  0 &  0 &  0 &  1 &  0 &  0 &  0 & -1 & -1 &  0 &  0 &  0 &  0 &  0 &  0 & -1 & -1 &  0 &  0 &  0 &  1 &  0 &  0 &  0 \\
 $ac$        &  0 &  0 &  0 &  0 &  1 &  0 &  0 &  0 &  0 &  0 & -1 & -1 &  0 &  0 &  0 & -1 &  0 & -1 &  0 &  0 &  0 &  0 &  1 &  0 \\
 $ab$        &  0 &  0 &  0 &  0 &  0 &  1 &  0 &  0 &  0 &  0 &  0 &  0 &  0 & -1 & -1 &  0 & -1 & -1 &  0 &  0 &  0 &  0 &  0 &  1 \\
 $d$         &  0 &  0 &  0 &  0 &  0 &  0 &  1 &  1 &  1 &  0 &  0 &  0 &  0 &  0 &  0 &  0 &  0 &  0 & -1 & -1 & -1 &  0 &  0 &  0 \\
 $c$         &  0 &  0 &  0 &  0 &  0 &  0 &  0 &  0 &  0 &  1 &  1 &  1 &  0 &  0 &  0 &  0 &  0 &  0 & -1 &  0 &  0 & -1 & -1 &  0 \\
 $b$         &  0 &  0 &  0 &  0 &  0 &  0 &  0 &  0 &  0 &  0 &  0 &  0 &  1 &  1 &  1 &  0 &  0 &  0 &  0 & -1 &  0 & -1 &  0 & -1 \\
 $a$         &  0 &  0 &  0 &  0 &  0 &  0 &  0 &  0 &  0 &  0 &  0 &  0 &  0 &  0 &  0 &  1 &  1 &  1 &  0 &  0 & -1 &  0 & -1 & -1 \\
 $\emptyset$ &  0 &  0 &  0 &  0 &  0 &  0 &  0 &  0 &  0 &  0 &  0 &  0 &  0 &  0 &  0 &  0 &  0 &  0 &  1 &  1 &  1 &  1 &  1 &  1 \\
\end{tabular}
}
\end{table}

It is instructive to write out all elementary imsets in a $2^{|N|}\times
|\cE(N)|$ matrix $\cfg_N$, 
where each elementary imset is a column of $\cfg_N$.
For $|N|=4$, $N=abcd$, the matrix $\cfg_N$ 
is written in Table \ref{tab:configuration}.
In this paper 
we  consider this matrix as the {\it configuration} defining a toric ideal.
Let $f(S)=\frac{1}{2} |S|^2$, $S\subseteq N$.  It is easily seen
that this $f$ is supermodular and 
\begin{equation}
\label{eq:homogeneous-configuration}
\sum_{S\subseteq N} f(S) u(S)=1, \quad   \forall u\in \cE(N) .
\end{equation}
Hence $\cfg_N$ defines a homogeneous toric ideal 
(Chapter 4 of Sturmfels\cite{sturmfels1996}).
Note that since $|S|$ is modular, any $f(S)=\frac{1}{2}|S|^2 + c_1 |S| + c_2$,
$c_1, c_2\in \RR$, satisfies 
\eqref{eq:homogeneous-configuration}.  In particular $f(S)=\frac{1}{2}|S|(|S|-1)$
is the standard supermodular function satisfying \eqref{eq:homogeneous-configuration}.

In Table \ref{tab:configuration} rows are ordered (from bottom to top)
by the cardinality of the set and by reverse lexicographic order among sets of
the same cardinality.  We call this order the graded reverse lexicographic
order of $\cP(N)$.  This order can be generalized to elementary imsets
as follows.  
In an elementary imset $u_{\tri{a}{b}{C}}$ we always order $a,b$ as $a<b$. 
Then we define
\begin{equation}
\label{eq:grlex-for-elementary}
u_{\tri{a}{b}{C}}\precw u_{\tri{a'}{b'}{C'}} \quad
\text{if}\quad
\begin{cases}
C<C', & \\
\text{or }  C=C' \text{ and } b<b',&\\
\text{or }  C=C' \text{ and } b=b' \text{ and } a<a' . &
\end{cases}
\end{equation}
In particular, strict inclusion $C\subset C'$ implies $C<C'$.
For $N=\{1,\dots,|N|\}$, 
the maximum element of $\cE(N)$ with respect to $\precw$ is $u_{\tri{1}{2}{R}}$, 
where $R=\{3,\ldots,|N|\}$.
We call this order the graded reverse lexicographic order of $\cE(N)$.
This order will be used in the proof of 
Theorem \ref{thm:kernel-lattice-basis} 
in Section \ref{sec:small-linear-relations}.
The columns in Table \ref{tab:configuration} are ordered 
(from right to left) according to the graded reverse lexicographic order
for elementary imsets.

\begin{figure}[htbp]
\setlength{\unitlength}{1.1mm}
\begin{center}
\begin{picture}(45,38)(0,7)
\put(20,10){\line(1,1){10}}
\put(20,10){\line(-1,1){10}}
\put(10,20){\line(1,1){10}}
\put(30,20){\line(-1,1){10}}
\put(30,20){\line(1,1){10}}
\put(40,30){\line(-1,1){10}}
\put(20,30){\line(1,1){10}}
\put(19,7){$C$}
\put(31,18){$BC$}
\put(40,29){$BCD$}
\put(3,19){$AC$}
\put(10,29){$ABC$}
\put(24,41){$ABCD$}
\put(19,12){$1$}
\put(11,19){$-1$}
\put(24,19){$-1$}
\put(19,26){$1$}
\put(29,22){$1$}
\put(21,29){$-1$}
\put(34,29){$-1$}
\put(29,36){$1$}
\end{picture}
\caption{Sum of two semi-elementary imsets}
\label{fig:two-quadrangles}
\end{center}
\end{figure}
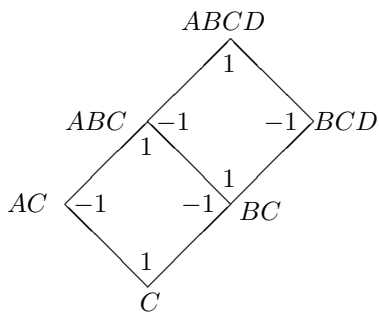

The following identity
among semi-elementary imsets is of basic importance.
\begin{equation}
\label{eq:sg-imset}
u_{\tri{A}{BD}{C}}=u_{\tri{A}{B}{C}}  + u_{\tri{A}{D}{BC}}. 
\end{equation}
This identity can be directly verified as
\begin{align*}
u_{\tri{A}{B}{C}}  + u_{\tri{A}{D}{BC}} &=
(\delta_{ABC}+\delta_C - \delta_{AC} - \delta_{BC})\\
& \qquad\quad + (\delta_{ABCD} + \delta_{BC} - \delta_{ABC} - \delta_{BCD})\\
&=\delta_{ABCD}  + \delta_C - \delta_{AC} - \delta_{BCD}.
\end{align*}
It is conveniently depicted in Fig.\ \ref{fig:two-quadrangles}.

From  \eqref{eq:sg-imset} we can split $A$ or $B$ in $u_{\tri{A}{B}{C}}$ into
smaller subsets.  If we repeat this splitting, 
every semi-elementary imset can be
written as a non-negative integer combination of
elementary imsets
\begin{equation}
\label{eq:semi-elementary-by-elementary}
u_{\tri{A}{B}{C}}=\sum_{v\in \cE(N)} k_v \cdot v,   \qquad \text{where}\ k_v \in \ZZ^+ .
\end{equation}
In Kashimura et al.\cite{kashimura-etal-2by2} we gave a detailed study of 
the set of all possible non-negative integer combinations of elementary imsets
which are equal to a semi-elementary imset.

We consider $\RR^{2^{|N|}}$ as equipped with the
standard inner product  $\langle \cdot, \cdot \rangle$.  Then
the inner product of $f: \cP(N)\rightarrow \RR^{2^{|N|}}$
and $u_{\tri{A}{B}{C}}$ is written as
\[
\langle f, u_{\tri{A}{B}{C}}\rangle = f(ABC)-f(C)+f(AC)+f(BC) .
\]
This inner product was already considered in \eqref{eq:homogeneous-configuration}.
Let $E=AC, F=BC, E\cap F=C, EF=ABC$ in \eqref{eq:super-modular-function}.
We see that $f$ is supermodular if and only if
\begin{equation}
\label{eq:dual-relation}
\langle f, u_{\tri{A}{B}{C}}\rangle \ge 0, \qquad \forall \tri{A}{B}{C}
\in \cT(N).
\end{equation}
Hence the set of semi-elementary imsets 
$\{ u_{\tri{A}{B}{C}}\mid \forall \tri{A}{B}{C} \in \cT(N)\}$
give the $H$-representation of the supermodular cone $\cK(N)$.
By definition, the convex cone 
\begin{align*}
\cK^*(N) &=\convcone \{ u_{\tri{A}{B}{C}} \mid \tri{A}{B}{C}\in \cT(N)\}\\
 &= \{ \sum_{\tri{A}{B}{C}\in \cT(N)} k_{\tri{A}{B}{C}} \cdot u_{\tri{A}{B}{C}}  
      \mid k_{\tri{A}{B}{C}} \ge 0 \}
\end{align*}
generated by the semi-elementary imsets is the cone dual to the
supermodular cone $\cK(N)$.  We call 
$\cK^*(N)$ the {\em imset cone}.  

In the $H$-representation of $\cK(N)$, not all of the hyperplanes 
determined by $u_{\tri{A}{B}{C}}$ are effective.
In fact by \eqref{eq:semi-elementary-by-elementary},
if $\langle f, v \rangle \ge 0$ for all elementary $v\in \cE(N)$, then
$\langle f, u_{\tri{A}{B}{C}} \rangle \ge 0$ for every semi-elementary
imset.  Hence $f$ is supermodular if 
$\langle f,v \rangle \ge 0$,
$\forall v\in \cE(N)$. 
It also
follows that $\cK^*(N)$ is generated by  the elementary imsets:
\begin{equation}
\label{eq:dual-cone}
\cK^*(N)=\convcone  \cE(N).
\end{equation}
$\cK^*(N)$ is a pointed cone by \eqref{eq:homogeneous-configuration}.
Note that \eqref{eq:dual-cone}
does not yet imply that the every elementary imset is indeed an
effective  hyperplane defining $\cK(N)$, or equivalently every elementary
imset is an extreme ray of $\cK^*(N)$.  
Note that an extreme ray of a polyhedral convex cone is a half line 
and when we say ``$u$ is an extreme ray of a cone'',
we actually mean that $\{ cu \mid c \ge 0\}$ is an extreme ray of the cone.
At this point we establish the following lemma.
\begin{lemma}
\label{lem:elementary-imset-extreme}
Each elementary imset is an extreme ray of the imset cone $\cK^*(N)$.
\end{lemma}
There are many ways to prove this lemma.  We give a somewhat involved 
argument, which will be used frequently in Section
\ref{sec:small-linear-relations}.
\begin{proof}
Let $u=u_{\tri{a}{b}{C}}\in \cE(N)$ be an elementary imset.
Suppose that $u$ is written as a non-negative  combination of 
elements of $\cK^*(N)$.  By \eqref{eq:dual-cone} 
\begin{equation}
\label{eq:level-proof-1}
u_{\tri{a}{b}{C}}= \sum_{v\in \cE(N)} k_v \cdot v, \qquad k_v \ge 0, \ \forall v\in \cE(N) .
\end{equation}
We need to show that $k_v=0$ for all $v\neq u_{\tri{a}{b}{C}}$.
Suppose that for $k_v > 0$ for some $v=u_{\tri{a'}{b'}{C'}}\neq u_{\tri{a}{b}{C}}$.
Among such $v$'s, choose $v_{\max}=u_{\tri{a'}{b'}{C'}}$ 
such that $a'b'C'$ is the largest according to the
graded reverse lexicographic order of $\cP(N)$.  
If $a'b'C' > abC$, then on the left-hand side
of \eqref{eq:level-proof-1}
$u_{\tri{a}{b}{C}}(a'b'C')=0$, 
but on the right-hand side of \eqref{eq:level-proof-1}
$\sum_{v\in \cE(N)} k_v \cdot v (a'b'C')>0$. This is a contradiction.
Hence if $k_v >0$, then $a'b'C' \le abC$.  Similarly,
considering the graded reverse lexicographic order of $C$, we have
$C' \ge C$, if $k_v > 0$ for $v=u_{\tri{a'}{b'}{C'}}$.
In particular for $v=u_{\tri{a'}{b'}{C'}}$ with $k_v > 0$ we have
$|C|\le |C'|$, $|a'b'C'|\le |abC|$.  But this implies  $|C|=|C'|$.
Then for any $S$ such that  $|S|=|C|$ we have
$\sum_{v\in \cE(N)} k_v \cdot v(S) \ge 0$.  Since $u_{\tri{a}{b}{C}}(S)=0$ for
$S\neq C$,  we have $k_v=0$ for $v=u_{\tri{a'}{b'}{S}}$, $S\neq C$, $|S|=|C|$. 
Similarly $k_v=0$ for $v=u_{\tri{a'}{b'}{C'}}$, $a'b'C'\neq abC$, $|C'|=|C|$.
Therefore $k_v >0$ only if $C'=C$ and $a'b'C'=abC$.  However this implies
$v=u_{\tri{a}{b}{C}}$.
\end{proof}

Let $\cL^*(N)\subseteq \RR^{2^{|N|}}$ denote the linear subspace 
spanned by the semi-elementary imsets. $\cL^*(N)$ is the orthogonal
complement of $\cL(N)$ and its dimension is
\[
\dim \cL^*(N)=2^{|N|} - |N|-1.
\]
$\cK^*(N) \subset \cL^*(N)$ and
the dimension of the relative interior of $\cK^*(N)$ is the
same as $\dim \cL^*(N)$.  By  \eqref{eq:modular-basis}, 
\[
u\in \cL^*(N) \quad \Leftrightarrow \quad
\sum_{E\in \cP(N)} u(E)=0, \ \  
\sum_{e\in E\in \cP(N)}u(E)=0, \ \forall e\in N.
\]

Now let
\begin{align*}
\cL_\ZZ^* &= \ZZ \{u_{\tri{A}{B}{C}}\mid \tri{A}{B}{C}\in \cT(N)\}\\
&=\{ \sum_{\tri{A}{B}{C}\in \cT(N)} k_{\tri{A}{B}{C}} \cdot u_{\tri{A}{B}{C}}
\mid k_{\tri{A}{B}{C}} \in \ZZ\}  \ \subseteq \ \ZZ^{2^{|N|}}
\end{align*}
denote the submodule of $\ZZ^{2^{|N|}}$
generated by the semi-elementary imsets.
Again by \eqref{eq:semi-elementary-by-elementary},
$\cL_\ZZ^*$ is 
generated by the elementary imsets:
\begin{equation}
\label{eq:integer-lattice1}
\cL_\ZZ^*=\ZZ \cE(N)=\{ \sum_{v\in \cE(N)} k_v \cdot v \mid k_v\in \ZZ\}.
\end{equation}
$\cL_\ZZ^*$ coincides with the set of integer points in $\cL^*(N)$:
\begin{equation}
\label{eq:integer-lattice2}
\cL_\ZZ^* = \cL^*(N) \cap \ZZ^{2^{|N|}}.
\end{equation}
This can be seen as follows.  For each $C\subseteq N$, $|C|\le |N|-2$,
choose an elementary imset $u_{\tri{a}{b}{C}}$. 
Consider a submatrix of the configuration consisting these elementary
imsets, which is of size $2^{|N|}\times (2^{|N|}-|N|-1)$. 
If $C$'s are ordered by the graded reverse lexicographic order 
as in Table \ref{tab:configuration}, 
then the lower $(2^{|N|}-|N|-1)\times (2^{|N|}-|N|-1)$ part of 
the matrix 
is an upper triangular matrix with 1 on the diagonal, as seen in Table
\ref{tab:configuration}. 
Therefore it is unimodular (= has integral inversion). 
Then the projection of $\cL_\ZZ^*$ onto these 
$2^{|N|}-|N|-1$ coordinates coincides with 
$\ZZ^{2^{|N|}-|N|-1}$.  Then \eqref{eq:integer-lattice2} holds, because
other elements are determined by orthogonality to $\cL(N)$.

Let $\QQ^+ \cE(N)$ denote the set of non-negative rational combinations
of elementary imsets.  
An integer point in $\QQ^+ \cE(N)$ 
is  called a {\em structural imset}. Let  $\cS(N)$ denote
the set of structural imsets.  Since $\cK(N)$ and $\cK(N)^*$ are 
rational polyhedral cones,  we have 
\[
\cS(N)= \cK^*(N)\cap \ZZ^{2^{|N|}}.
\]
By \eqref{eq:integer-lattice2}, $\cS(N)$ can also be written as $\cS(N)=\cK^*(N)\cap \cL_\ZZ^*$.

A non-negative integer combination of elementary imsets is called a
{\em combinatorial imset} and    
\[
\cC(N) = \ZZ^+ \cE(N)
= \{ \sum_{v\in \cE(N)} k_v \cdot v \mid k_v \in \ZZ^+\}
\]
denotes the set of combinatorial imsets.
$\cC(N)$ is the semigroup  generated by $\cE(N)$.
Clearly $\cC(N)\subseteq \cS(N)$ and it is known that for $|N|\ge 5$
the inclusion is strict (Hemmecke et al.\cite{hemmecke-etal-2008}), i.e.\ for $|N|\ge 5$ the
semigroup $\cC(N)$ is not normal.

\section{Multiinformation and derivation of conditional independence using imsets}
\label{sec:multiinformation}

Let $P,\mu$ be two probability measures on a sample space ${\cal X}$, such 
that $P$ is absolutely continuous with respect to $\mu$.
The relative entropy $H(P|\mu)$ of $P$ with respect to $\mu$ is defined
by
\[
H(P|\mu)=\int_{\cal X} \ln \frac{{\rm d}P}{{\rm d}\mu}(x)  \; {\rm d}P(x).
\]
$H(P|\mu)\ge 0$ and $H(P|\mu)=0$ if and only if $P=\mu$.

Consider a joint probability distribution $P$ of  the variables in $N$.
As usual, 
$\indtri{A}{B}{C}{P}$ denotes the conditional independence statement of
variables in $A$ and in $B$ given the variables in $C$ under $P$.
The set of conditional independences under $P$ is denoted by
\[
\cM_P = \bra{ \tri{A}{B}{C} \in \mathcal{T}(N) \mid \indtri{A}{B}{C}{P}}. 
\]
We call $\cM_P$ the conditional independence model of $P$.

For $A\subseteq N$, $P^A$ denotes the marginal distribution of the
variables in $A$.   
Let $\prod_{i\in A} P^{\{i\}}$ denote the
product of one-dimensional marginal distributions in $A$.  Then
the multiinformation function $m_P: \cP(N) \rightarrow [0,\infty]$ is defined
by
\begin{equation}
\label{eq:multiinformation}
m_P(S)= H(P^S \mid \prod_{i\in S} P^{\{i\}}),   \ S\neq\emptyset,
\end{equation}
and $m_P(\emptyset)=0$.  Throughout this paper we only consider
$P$ such that $m_P(S)$ is finite for every $S\subseteq N$.
The basic fact on the multiinformation is that $m_P$ is supermodular. 
Furthermore
the following equivalence between the conditional independence 
$\indtri{A}{B}{C}{P}$
and the local modularity of $m_P$ at  $\tri{A}{B}{C}$ holds:
\begin{align}
\langle m_P, u_{\tri{A}{B}{C}}\rangle
&\ge 0,   \ \ \ \forall \tri{A}{B}{C}\in \cT(N)\nonumber\\
                               &=0    \ \ \Leftrightarrow \ \ 
\indtri{A}{B}{C}{P}.
\label{eq:multiinformation-basic}
\end{align}
This equivalence is the basis for manipulating conditional independence statements
in terms of imsets. 
Note that  $m_P$ is
standardized as a supermodular function
since $m_P(S)=0$ for $|S|\le 1$.


Traditionally, the implications among conditional independence statements under a probability 
measure $P$ have been
studied in terms of the following semi-graphoid axioms.  
In the axioms $A,B,C,D\subseteq N$ are disjoint.

\bigskip
\setlength{\tabcolsep}{2pt}
\begin{tabular}{lll}
1.&triviality & $\indtri{A}{\emptyset}{C}{P}$\\
2.&symmetry   & $\indtri{A}{B}{C}{P} \ \Rightarrow \ \indtri{B}{A}{C}{P}$ \\
3.&decomposition & $\indtri{A}{BD}{C}{P} \ \Rightarrow\ \indtri{A}{B}{C}{P}$\\
4.& weak union &  $\indtri{A}{BD}{C}{P} \ \Rightarrow\ \indtri{A}{D}{BC}{P}$\\
5.& contraction & $\indtri{A}{D}{BC}{P} \ \text{and}\ \indtri{A}{B}{C}{P} \Rightarrow 
\indtri{A}{BD}{C}{P}$
\end{tabular}

\bigskip
Note that decomposition, weak union and contraction can be combined into the
following single equivalence:
\[
\indtri{A}{BD}{C}{P} \ \ \Leftrightarrow\ \ 
\indtri{A}{D}{BC}{P} \ \text{and}\ \indtri{A}{B}{C}{P} .
\]
This equivalence can be proved very easily by imsets.
Take the inner product of 
\eqref{eq:sg-imset} with the multiinformation $m_P$.  Then
\[
\langle m_P, u_{\tri{A}{BD}{C}} \rangle
= \langle m_P, u_{\tri{A}{B}{C}}  \rangle
 + \langle m_P, u_{\tri{A}{D}{BC}} \rangle.
\]
Since every term is non-negative, we have
\begin{align*}
\indtri{A}{BD}{C}{P}  &\Leftrightarrow
\langle m_P, u_{\tri{A}{BD}{C}} \rangle = 0 \\
&\Leftrightarrow
\langle m_P, u_{\tri{A}{B}{C}} \rangle 
= \langle m_P, u_{\tri{A}{D}{BC}}\rangle =0 \\
&\Leftrightarrow
\indtri{A}{D}{BC}{P} \ \text{and}\ \indtri{A}{B}{C}{P} .
\end{align*}
In this way, manipulation of ``rules'' such as the semi-graphoid axioms
is translated to linear algebraic operations in terms of imsets.
This is a very important advantage of the method of imsets and the multiinformation.

As another example, consider the following identity.
\begin{align}
u_{\tri{a}{b}{c}} + u_{\tri{a}{c}{d}} + u_{\tri{a}{d}{b}}
&= (\delta_{abc} + \delta_c - \delta_{ac} - \delta_{bc})
 + (\delta_{acd} + \delta_d - \delta_{ad} - \delta_{cd}) \nonumber\\
& \qquad + (\delta_{abd} + \delta_b - \delta_{ab} - \delta_{bd}) \nonumber\\
&=(\delta_{abc}+\delta_{acd}+\delta_{abd}) + (\delta_b + \delta_c + \delta_d)
\nonumber\\
&\qquad - (\delta_{ab} + \delta_{ac}+\delta_{ad}) - (\delta_{bc} + \delta_{bd} + \delta_{cd})\nonumber \\
&=u_{\tri{a}{c}{b}}+u_{\tri{a}{d}{c}}+u_{\tri{a}{b}{d}}.
\label{eq:3:3-relation}
\end{align}
By the method of imsets we then have
\begin{equation}
\label{eq:3-3-equivalence}
\indtr{a}{b}{c}, 
\indtr{a}{c}{d}, 
\indtr{a}{d}{b}
\ \Leftrightarrow\ 
\indtr{a}{c}{b}, 
\indtr{a}{d}{c}, 
\indtr{a}{b}{d},
\end{equation}
where for simplicity we wrote $\indtr{a}{b}{c}
$ instead of $\indtri{a}{b}{c}{P}$.
It is evident that none of decomposition, weak union and contraction
can be applied to the left-hand side nor to the right-hand side of this
equivalence.  Therefore this equivalence can not be derived from semi-graphoid
axioms.

Here it is interesting to note that \eqref{eq:3:3-relation} 
can be derived by linear algebraic operations from 
\eqref{eq:sg-imset}, which itself corresponds to the semi-graphoid axiom.
By applying \eqref{eq:sg-imset} twice we can write
\begin{align}
3 \cdot u_{\tri{a}{bcd}{\emptyset}}
&=(u_{\tri{a}{b}{\emptyset}} + u_{\tri{a}{c}{b}} + u_{\tri{a}{d}{bc}})\nonumber\\ 
&\quad + (u_{\tri{a}{c}{\emptyset}} + u_{\tri{a}{d}{c}}  + u_{\tri{a}{b}{cd}})\nonumber\\
&\quad + (u_{\tri{a}{d}{\emptyset}} + u_{\tri{a}{b}{d}} + u_{\tri{a}{c}{bd}})
\label{eq:3-3-1}
\end{align}
as well as
\begin{align}
3 \cdot u_{\tri{a}{bcd}{\emptyset}}
&=(u_{\tri{a}{b}{\emptyset}} + u_{\tri{a}{d}{b}} + u_{\tri{a}{c}{bd}}) \nonumber\\ 
&\quad + 
(u_{\tri{a}{c}{\emptyset}} + u_{\tri{a}{b}{c}} + u_{\tri{a}{d}{bc}}) \nonumber\\
&\quad + (u_{\tri{a}{d}{\emptyset}} + u_{\tri{a}{c}{d}} + u_{\tri{a}{b}{cd}}).
\label{eq:3-3-2}
\end{align}
Equality of right-hand sides of \eqref{eq:3-3-1} and \eqref{eq:3-3-2}
immediately gives  \eqref{eq:3:3-relation}.

Let $u\in \cS(N)$ be a structural imset.  A conditional independence with respect to $u$
is defined as follows.
We say that $A$ and $B$ are conditionally independent given $C$ with respect to $u$,  
and denote it by $\indtri{A}{B}{C}{u}$, if
\begin{equation}
\label{eq:ci-under-u}
\exists k \in \NN, \ k \cdot u - u_{\tri{A}{B}{C}} \in \mathcal{S}(N). 
\end{equation}
The set of conditional independences  induced by $u$ is denoted by
\begin{equation}
\label{eq:indep-u}
\mathcal{M}_u = \bra{ \tri{A}{B}{C} \in \mathcal{T}(N) \mid \indtri{A}{B}{C}{u}}. 
\end{equation}
We call $\mathcal{M}_u$ the conditional independence model of $u$.
%
The importance of this definition lies in the following completeness theorem for imsets.
\begin{theorem}[Theorem 5.2 of Studen\'y\cite{stu2005}]
\label{thm:complete}
Let $P$ be a probability measure over $N$ with finite multiinformation.
Then there exists $u\in \cS(N)$ such that $\cM_P = \cM_u$.
\end{theorem}

\begin{proof}
Define $u=\sum_{\tri{A}{B}{C}\in \cM_P} u_{\tri{A}{B}{C}}$.
Then clearly  $\cM_P \subseteq \cM_u$.  Now for  any $\tri{A}{B}{C}\in \cM_u$
choose $k$ such that 
\eqref{eq:ci-under-u} holds.
Then 
\[
0=\langle m_P, k\cdot u\rangle
=\langle m_P, k\cdot u - u_{\tri{A}{B}{C}} \rangle + \langle m_P, u_{\tri{A}{B}{C}} \rangle.
\]
Since two terms on the right-hand side are non-negative, we have
$\langle m_P, u_{\tri{A}{B}{C}} \rangle=0$ and hence  $\cM_u \subseteq \cM_P$.
\end{proof}

We have reproduced this proof from Section 5 of {Studen\'y\cite{stu2005}}, because
this argument is instructive.  In this proof, to show 
$\langle m_P, u_{\tri{A}{B}{C}} \rangle=0$, we have added an extra term
$k\cdot u - u_{\tri{A}{B}{C}}$ to $u_{\tri{A}{B}{C}}$.  In this sense, this 
proof is similar to  the argument concerning  \eqref{eq:3-3-1}
and \eqref{eq:3-3-2} above for showing \eqref{eq:3-3-equivalence}.

Next we show that the conditional independence model
$\cM_u$ depends only on the 
face of $\cK^{*}(N) $ 
which contains $u$ as a relatively interior point.
Although this fact is discussed in a series of papers by Studen\'y\cite{stu1994}
and Chapter 5 of {Studen\'y\cite{stu2005}}, 
we make this point clear in the following lemma.

\begin{lemma} \label{lem:cone}
Let $a_1, \dots, a_n$ be nonzero vectors in $\RR^m$.
Let $C=\convcone\{a_1, \dots, a_n\}$ be the convex cone generated by these vectors
and suppose that $C$ is pointed and $a_1,\dots,a_n$ are extreme rays of $C$.
For a face $F$ of $C$, let $F^\circ$ denote its relative interior.
Put  $I_F = \bra{ i \mid a_i \in F}$.
Then, $b \in C$ belongs to $F^\circ$ if and only if the following two conditions hold: 
\begin{enumerate}
\item For every $i \in I_F$, there exists $k_i \in \NN$ such that
$  k_i \cdot b - a_i \in C$, 
\item For every $i \not\in I_F$ and every  $k \in \NN$,
$ k \cdot b - a_i \not\in C.$
\end{enumerate}
\end{lemma}

\begin{proof}
$b\in C$ can be written as
$b = \sum_{i=1}^n c_i a_i$, $c_i \ge 0$, $i= 1, \dots, n.$ 
Suppose that $b$ satisfies the above two conditions. 
Then for $i \in I_F$, 
by condition 1, for some $k_i \in \NN$, $k_i b = a_i + \sum_{j=1}^n \tilde c_{j}^{i} a_j$, $\tilde c_{j}^{i} \ge 0$, or
\[
b = \frac{1}{k_i} a_i + \sum_{j=1}^n (\tilde c_{j}^{i}/k_i) a_j.
\]
Taking the average of the right-hand side over  $i\in I_F$, we see that
we can take $c_i > 0$ for $i\in I_F$ in $b = \sum_{i=1}^n c_i a_i$.
On the other hand,  for $i \not\in I_F$, 
condition 2 implies $k \cdot c_i - 1 \leq 0$ for any $k \in \NN$, hence  $c_i = 0$.
Therefore if 
$b$ satisfies two conditions then 
$b$ is written as
\[ b =\sum_{i \in I_F} c_i a_i, \quad c_i > 0, i \in I_F. \]
Since the coefficients are positive, we have $b \in F^\circ$.

Conversely let $b \in F^\circ$. For $i \in I_F$ and sufficiently large $k \in \NN$, 
let 
$a = (1/k) a_i$.
Then we have $b - a \in F$, which is equivalent to $k \cdot b - a_i \in F \subseteq C$.
Let $v \in \RR^m$ be a normal vector of a supporting hyperplane of $F$ such that 
$ \langle v, x \rangle \ge 0$, $\forall x \in C,$ and  
for $b \in F$ and any $a_i \not\in F$ the inner product with respect to $v$ satisfy  
\begin{align*}
        \langle v, b \rangle = 0, \quad \langle v, a_i \rangle > 0.
\end{align*}
Then for any $k \in \NN$, we have  
\begin{align*}
        \langle v, k \cdot b - a_i \rangle < 0,
\end{align*}
which means that $k \cdot b - a_i \not\in C$.
\end{proof}

Under the assumptions of Lemma \ref{lem:cone},
let $b\in C$ and let $F=F_b$ denote the face of $C$ such that
$b \in F^\circ$. Let  
\[
\cE_b = \{ a_i \mid k\cdot b - a_i \in C, \  \exists k\in \NN\}.
\]
Then $\cE_b$ is the set of extreme rays of
$F_b$ by Lemma \ref{lem:cone}.

{}From this fact
we can prove that conditional independence structures induced by imsets
depend only on faces of the imset cone $\cK^*(N)$ and not on each imset. 
\begin{proposition}
\label{prop:face}
Let $u, u' \in \mathcal{S}(N)$.
Then $\mathcal{M}_u = \mathcal{M}_{u'}$ if and only if $u, u'$ belong to the relative interior  of the same face of $\cK^*(N)$. 
\end{proposition}

\begin{proof} For $v\in \cS(N)$,  let $F_v$ denote the face of $\cK^*(N)$ such that 
$v$ belongs to the relative interior  $F_v^\circ$ of $F_v$ and
let $\mathcal{E}_v(N)$ denote 
the set of elementary imset $u_{\tri{a}{b}{C}}$ such that 
$\exists k \in \NN, \ k \cdot v - u_{\tri{a}{b}{C}} \in \mathcal{S}(N)$,
that is, 
\[ 
\mathcal{E}_v(N) = \bra{ u_{\tri{a}{b}{C}} \mid \tri{a}{b}{C} \in \mathcal{M}_u}.
\]
As note above, 
$\mathcal{E}_v(N)$ is the set of extreme rays of
$F_v$. Hence
\begin{align*}
\mathcal{M}_u = \mathcal{M}_{u'} & \Leftrightarrow  \  \mathcal{E}_u(N) = \mathcal{E}_{u'}(N)\\
& \Leftrightarrow \ F_u = F_{u'}.
\end{align*}
\end{proof}

By Theorem \ref{thm:complete} and Proposition \ref{prop:face}, each
conditional independence model $\cM_P$ corresponds to 
a face of $\cK^*(N)$.  
Note that, by mutual orthogonality, the face poset of $\cK(N)$ and that of 
$\cK^*(N)$ are isomorphic (when inclusion in $\cK(N)$ is reversed).
Hence each conditional independence model $\cM_P$ also corresponds to 
a face of $\cK(N)$.  

Also note that 
we can regard the set of conditional independence models $\cM_P$ 
as a poset with respect to the inclusion relation
$\cM_P \subseteq \cM_{P'}$.
By Proposition \ref{prop:face} this poset is  a sub-poset
of the face poset of $\cK^*(N)$.  Unfortunately it is known that
for $|N|\ge 4$, there exist faces of $\cK^*(N)$ which do not
correspond to any $\cM_P$, i.e.\ for $|N|\ge 4$, 
the set of conditional independence models is a proper sub-poset
of the face poset of $\cK^*(N)$.

\begin{figure}[htbp]
\caption{An extreme ray of $\cK(N)$ not corresponding to multiinformation}
\label{fig:extreme-ray}
\setlength{\unitlength}{0.5mm}
\begin{center}
\begin{picture}(100,120)(-50,-10)
\put(0,-5){\oval(24,16)[c]}
\put(-12,-6){\line(1,0){24}}
\put(-5,-15){\makebox(10,10)[c]{$\emptyset$}}
\put(-5,-5){\makebox(10,6)[c]{0}}
\put(0,3){\line(-5,1){45}}
\put(0,3){\line(-3,2){14}}
\put(0,3){\line(3,2){14}}
\put(0,3){\line(5,1){45}}
\put(-45,20){\oval(24,16)[c]}
\put(-57,19){\line(1,0){24}}
\put(-50,11){\makebox(10,10)[c]{$a$}}
\put(-50,20){\makebox(10,6)[c]{0}}
\put(-15,20){\oval(24,16)[c]}
\put(-27,19){\line(1,0){24}}
\put(-20,11){\makebox(10,10)[c]{$b$}}
\put(-20,20){\makebox(10,6)[c]{0}}
\put(15,20){\oval(24,16)[c]}
\put(3,19){\line(1,0){24}}
\put(10,11){\makebox(10,10)[c]{$c$}}
\put(10,20){\makebox(10,6)[c]{0}}
\put(45,20){\oval(24,16)[c]}
\put(33,19){\line(1,0){24}}
\put(40,11){\makebox(10,10)[c]{$d$}}
\put(40,20){\makebox(10,6)[c]{0}}
\put(-75,45){\oval(24,16)[c]}
\put(-87,44){\line(1,0){24}}
\put(-80,35){\makebox(10,10)[c]{$ab$}}
\put(-80,45){\makebox(10,6)[c]{0}}
\put(-45,45){\oval(24,16)[c]}
\put(-57,44){\line(1,0){24}}
\put(-50,36){\makebox(10,10)[c]{$ac$}}
\put(-50,45){\makebox(10,6)[c]{1}}
\put(-15,45){\oval(24,16)[c]}
\put(-27,44){\line(1,0){24}}
\put(-20,36){\makebox(10,10)[c]{$ad$}}
\put(-20,45){\makebox(10,6)[c]{1}}
\put(15,45){\oval(24,16)[c]}
\put(3,44){\line(1,0){24}}
\put(10,36){\makebox(10,10)[c]{$bc$}}
\put(10,45){\makebox(10,6)[c]{1}}
\put(45,45){\oval(24,16)[c]}
\put(33,44){\line(1,0){24}}
\put(40,36){\makebox(10,10)[c]{$bd$}}
\put(40,45){\makebox(10,6)[c]{1}}
\put(75,45){\oval(24,16)[c]}
\put(63,44){\line(1,0){24}}
\put(70,36){\makebox(10,10)[c]{$cd$}}
\put(70,45){\makebox(10,6)[c]{1}}
%
%
\put(-73,37){\line(6,-1){56}}
\put(-74,37){\line(3,-1){27}}
\put(-73,53){\line(6,1){56}}
\put(-74,53){\line(3,1){27}}
\put(-43,37){\line(6,-1){56}}
\put(-45,37){\line(0,-1){9}}
\put(-43,53){\line(6,1){56}}
\put(-45,53){\line(0,1){9}}
\put(-15,37){\line(-3,-1){27}}
\put(-14,37){\line(6,-1){56}}
\put(-15,53){\line(0,1){9}}
\put(-14,53){\line(3,1){27}}
\put(14,37){\line(-3,-1){27}}
\put(15,37){\line(0,-1){9}}
\put(14,53){\line(-6,1){56}}
\put(15,53){\line(3,1){27}}
\put(43,37){\line(-6,-1){56}}
\put(45,37){\line(0,-1){9}}
\put(43,53){\line(-6,1){56}}
\put(45,53){\line(0,1){9}}
\put(73,37){\line(-6,-1){56}}
\put(74,37){\line(-3,-1){27}}
\put(73,53){\line(-6,1){56}}
\put(74,53){\line(-3,1){27}}
\put(-45,70){\oval(24,16)[c]}
\put(-57,69){\line(1,0){24}}
\put(-50,61){\makebox(10,10)[c]{$abc$}}
\put(-50,70){\makebox(10,6)[c]{2}}
\put(-15,70){\oval(24,16)[c]}
\put(-27,69){\line(1,0){24}}
\put(-20,61){\makebox(10,10)[c]{$abd$}}
\put(-20,70){\makebox(10,6)[c]{2}}
\put(15,70){\oval(24,16)[c]}
\put(3,69){\line(1,0){24}}
\put(10,61){\makebox(10,10)[c]{$acd$}}
\put(10,70){\makebox(10,6)[c]{2}}
\put(45,70){\oval(24,16)[c]}
\put(33,69){\line(1,0){24}}
\put(40,61){\makebox(10,10)[c]{$bcd$}}
\put(40,70){\makebox(10,6)[c]{2}}
\put(0,95){\oval(24,16)[c]}
\put(-12,94){\line(1,0){24}}
\put(-5,86){\makebox(10,10)[c]{$abcd$}}
\put(-5,95){\makebox(10,6)[c]{4}}
\put(0,87){\line(-5,-1){45}}
\put(0,87){\line(-3,-2){14}}
\put(0,87){\line(3,-2){14}}
\put(0,87){\line(5,-1){45}}
\end{picture}
\end{center}
\end{figure}
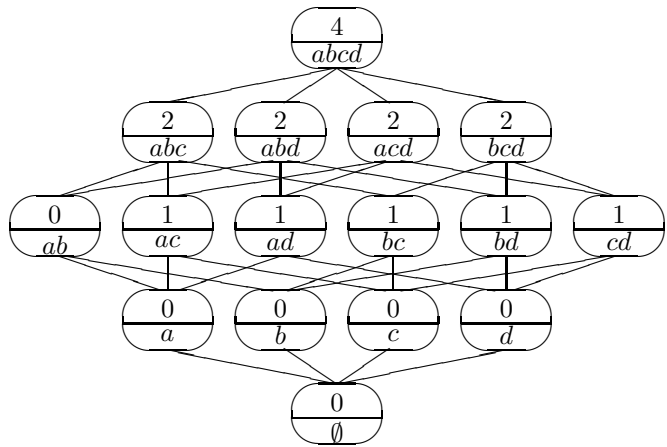

\begin{example}[Example 4.1 of Studen\'y\cite{stu2005}]
\label{ex:ex1}
For $N=abcd$ consider 
\begin{align*}
u &= u_{\tri{c}{d}{ab}} + u_{\tri{a}{b}{\emptyset}} 
+ u_{\tri{a}{b}{c}} + u_{\tri{a}{b}{d}}\\
&=\delta_N + 2\cdot \delta_{ab} -\delta_{ac} - \delta_{ad}
 - \delta_{bc} - \delta_{bd} -\delta_a - \delta_b + \delta_c + \delta_d 
 + \delta_\emptyset.
\end{align*}
Then $\tri{c}{d}{ab}, \tri{a}{b}{\emptyset}, \tri{a}{b}{c}, \tri{a}{b}{d}\in 
\cM_u$.
Furthermore let
\begin{equation}
\label{eq:reino-extreme-ray}
m = 4\cdot \delta_N + 2 \cdot \sum_{S\subseteq N: |S|=3}\delta_S
+\sum_{S\subseteq N: |S|=2, S\neq ab} \delta_S .
\end{equation}
(see Fig.\ \ref{fig:extreme-ray}.)
Then it can be checked that $\langle m, v \rangle \ge 0$ for
every $v\in \cE(N)$. Hence $m$ is supermodular.
Furthermore $\langle m, u\rangle  = 0$, $\langle m, u_{\tri{a}{b}{cd}}\rangle
= 1$. Consequently for any $k\in \NN$,  
$\langle m, k \cdot u - u_{\tri{a}{b}{cd}}\rangle = -1$ and hence
$k\cdot u -  u_{\tri{a}{b}{cd}} \not\in \cS(N)$ 
and $\tri{a}{b}{cd}\not\in \cM_u$.
On the other hand it can be shown that
for any probability measure $P$ with finite multiinformation
\[
\tri{c}{d}{ab}, \tri{a}{b}{\emptyset}, \tri{a}{b}{c}, \tri{a}{b}{d} \in \cM_P
\ \Rightarrow \ \tri{a}{b}{cd}\in \cM_P
\]
(cf. Corollary 2.1 of Studen\'y\cite{stu2005}). Therefore for the above $u$, there exists no $P$
such that $\cM_u =\cM_P$.
\end{example}
It can be shown that $m$ in  \eqref{eq:reino-extreme-ray} is an
extreme ray of the supermodular cone $\cK(N)$ for $|N|=4$ in the sense of the next section.

\section{Results on extreme rays of the supermodular cone}
\label{sec:extreme-rays}

In this section we prove that certain supermodular functions are
``extreme rays'' of the supermodular cone  $\cK(N)$.
Since extreme rays of $\cK(N)$ are hard to describe, even partial results on
the extreme rays are of interest.
As we discussed in Section \ref{sec:introduction}, 
since $\cK(N)$
contains the linear space $\cL(N)$ of modular functions,
in order to consider ``extreme rays'' of $\cK(N)$ we need
to identify two supermodular functions $f$ and $g$ if they differ
by a modular function. 
For this purpose it is
simplest to consider the cone $\cK_\ell(N)$ of standardized supermodular functions
in \eqref{eq:super-modular-standardization}:
\[
\cK_\ell(N)= \{ f\in \cK(N) \mid f(S)=0 \ \text{for}\ |S|\le 1\},
\]
which is a pointed polyhedral cone in ${\mathbb R}^{2^{|N|}-|N|-1}$.

Following the terminology of {Studen\'y\cite{stu2005}}, we call an extreme ray
of $\cK_\ell(N)$ a {\em skeletal} supermodular function.   
Furthermore by abusing the terminology, we call $f\in \cK(N)$ a skeletal
supermodular function (or an extreme ray of $\cK(N)$) 
if its standardization $\bar f$ in \eqref{eq:super-modular-standardization}
is a skeletal supermodular function.
Below we  often omit ``supermodular'' and
simply say that $f$ is a skeletal function, or  $f$ is skeletal.
Skeletal functions for $|N|=4$ are
depicted in the appendix of  Studen\'y et al.\cite{studeny-etal-five-variables}.
Various conditions on skeletal functions of this section are
meant to be useful for understanding these functions.
Except for the extreme ray 
$m$ in \eqref{eq:reino-extreme-ray}, 
other extreme rays for $|N|=4$ are covered by the results of this section.
However they are far from enough for understanding most of skeletal
functions for $|N|=5$ studied in Studen\'y et al\cite{studeny-etal-five-variables}.

The results of this section concern extending a supermodular function on ${\cal P}(N)$
to a larger base set $\tilde N \supset N$ and resemble results in 
Chapter 2 of  {Topkis\cite{topkis1998}}, although Topkis\cite{topkis1998}
does not consider skeletal functions.   Results on a skeletal function
as a maximum of a collection of modular functions are given in
Rosenm\"uller and  Weidner\cite{rosenmuller-weidner-1974}.

For proving that a 
supermodular function $f$ is skeletal
we proceed as in  Lemma \ref{lem:elementary-imset-extreme}.
Express $f$ as a sum of two supermodular functions
\begin{equation}
\label{eq:supermodular-sum}
f=f_1 + f_2.
\end{equation}
Here we can assume that $f, f_1, f_2$ are standardized, namely,
$0=f(S)=f_1(S)=f_2(S)$ for $|S|\le 1$. This follows from the fact that
the standardization \eqref{eq:super-modular-standardization} can be performed
to $f, f_1, f_2$ separately.
We need to show that $f_1$ and $f_2$ are proportional to $f$, i.e.,
\begin{equation}
\label{eq:supermodular-extreme-condition}
f_1=c_1  f, \ f_2 = c_2  f, \quad c_1, c_2 \ge 0,  \ c_1 + c_2 = 1.
\end{equation}

We first state two simple facts on skeletal  functions.
The first one concerns complementation of subsets of $N$.  
In Section 9.1.2 of Studen\'y\cite{stu2005} this is referred to as
a reflection.
Let $A^C$ denote the complement of $A\subseteq N$.
Let $f$ be a supermodular function.  Define $f^C$ by
\[
f^C(S)=f(S^C), \quad S\subseteq N.
\]
Taking the complement corresponds to looking at the configuration
$\cfg_N$ in Table \ref{tab:configuration} upside down.
The following lemma is trivial.
\begin{lemma}
\label{lem:complementation}
$f^C$ is  skeletal if and only if $f$ is 
skeletal.
\end{lemma}
The second one concerns a supermodular function $f(S)$ depending only on the
size $|S|$ of $S$.
\begin{lemma}
\label{lem:max-k}
Let $|N|\ge 2$ and $1\le k < |N|$.  $f(S)=\max(|S|-k,0)$ is skeletal.
\end{lemma}
\begin{proof}
 It is easily verified that $f$ is supermodular. 
 Let $f=f_1 + f_2$, where $f_1, f_2$ are standardized and supermodular. 
 Since $f_1, f_2$ are non-negative, $f_{1}=f_{2}=0$ holds for $|S|\ \leq k$. 
 This means \eqref{eq:supermodular-extreme-condition} holds for $|S|\ \leq k$. 
 Next we consider the case for $\abs{S} \geq k$. 
 If $k = \abs{N}-1$, then \eqref{eq:supermodular-extreme-condition} holds. 
 Hence assume $k \leq \abs{N}-2$. 
 Let $S_{k}$ be any subset of $N$ such that $\abs{S_{k}}=k$. 
  Let $\tilde{S} \supset S_{k}$. 
 Then $f(\tilde{S})=f(T \cup S_{k})$ where $T = \tilde{S} \setminus S_{k}$. 
 For $S_{k} \subseteq N$, define $\tilde{f}^{S_{k}}(T) = f(T \cup S_{k}) = \abs{T},\, T \subseteq N \setminus S_{k}$. 
 In the same way, define $\tilde{f}_{i}^{S_{k}}(T)=f_{i}(T \cup S_{k}),\, T \subseteq N \setminus S_{k}$, $i=1,2$. 
 Then $\tilde{f}^{S_{k}}$ and $\tilde{f}_{i}^{S_{k}}$ are modular function on $N \setminus S_{k}$. 
 By (\ref{eq:modular}), $\tilde{f}_{1}^{S_{k}}$ and $\tilde{f}_{2}^{S_{k}}$ can be written as 
 \begin{eqnarray}
  \tilde{f}_{1}^{S_{k}}(T) = \lambda_{\emptyset}^{S_{k}} + \sum_{a \in T}\lambda_{a}^{S_{k}}
   ,\quad
  \tilde{f}_{2}^{S_{k}}(T) = \mu_{\emptyset}^{S_{k}} + \sum_{a \in T}\mu^{S_{k}}. 
  \nonumber 
 \end{eqnarray}
 From $\tilde{f}^{S_{k}}(T) = \abs{T}$ and $f=f_1 + f_2$, we have 
 $\lambda_{\emptyset}^{S_{k}}=\mu_{\emptyset}^{S_{k}}=0$ and $\lambda_{\alpha}^{S_{k}}+\mu_{\alpha}^{S_{k}}=1$ for $\alpha  \in T$. 
 If $\lambda_{\alpha}^{S_{k}} = \lambda_{\beta}^{S_{k}}=\lambda^{S_{k}}$ and $\mu_{\alpha}^{S_{k}} = \mu_{\beta}^{S_{k}} = \mu^{S_{k}}$ 
 hold for all $\alpha, \beta \in N \setminus S_{k}$ where $\lambda^{S_{k}}, \mu^{S_{k}} \in \RR^{+}$, 
 then $\tilde{f}_{1}^{S_{k}}(T)=\lambda^{S_{k}}\abs{T}$ and $\tilde{f}_{2}^{S_{k}}(T)=\mu^{S_{k}}\abs{T}$ hold. 
 Furthermore we obtain the representations $\tilde{f}_{1}^{S_{k}}(T)=\lambda\abs{T}$ and $\tilde{f}_{2}^{S_{k}}(T)=\mu\abs{T}$ where $\lambda, \mu \in \RR^{+}$,
 because $\tilde{f}_{i}^{S_{k}}(\alpha')=\tilde{f}_{i}^{S_{k}'}(\alpha)$ holds if $\alpha' S_{k} = \alpha S_{k}'$ 
 where $S_{k}'$ is a subset of $N$ such that $\abs{S_{k}}'=k$ 
 and $\alpha \in S_{k},\, \alpha' \in S_{k}'$. 
 This means that \eqref{eq:supermodular-extreme-condition} holds for $\abs{S} \geq k$. 
 It is suffices to show that 
 $\lambda_{\alpha}^{S_{k}} = \lambda_{\beta}^{S_{k}}$ holds for all $\alpha, \beta \in N \setminus S_{k}$ such that $\alpha \neq \beta$, 
 because we can prove the equation $\mu_{\alpha}^{S_{k}} = \mu_{\beta}^{S_{k}}$ in the same manner. 
 Let $\alpha, \beta \in N \setminus S_{k}$ such that $\alpha \neq \beta$. 
 For $\gamma \in S_{k}$, let $S_{k}'=\{S\setminus \gamma\}\cup\alpha$ and $S_{k}''=\{S\setminus \gamma\}\cup\beta$. 
 First we consider $\tilde{f}_{1}^{S_{k}'}$. 
 From $\gamma S_{k}'=\alpha S_{k}$, we have $\tilde{f}_{1}^{S_{k}'}(\gamma)=\tilde{f}_{1}^{S_{k}}(\alpha)=\lambda_{\alpha}^{S_{k}}$. 
 Furthermore, from $\beta\gamma S_{k}'=\alpha\beta S_{k}$, 
 we have $\tilde{f}_{1}^{S_{k}'}(\beta\gamma)=\tilde{f}_{1}^{S_{k}}(\alpha\beta)=\lambda_{\alpha}^{S_{k}}+\lambda_{\beta}^{S_{k}}$. 
 Since $\tilde{f}_{1}^{S_{k}'}$ is modular on $N \setminus S_{k}'$, 
 we obtain 
 \begin{equation}
  \tilde{f}_{1}^{S_{k}'}(\beta)
   = \tilde{f}_{1}^{S_{k}'}(\beta\gamma) + \tilde{f}_{1}^{S_{k}'}(\emptyset) - \tilde{f}_{1}^{S_{k}'}(\gamma)
   =\lambda_{\beta}^{S_{k}}. 
   \label{eq:f1Skd_modular}
 \end{equation}
 Next we consider $\tilde{f}_{1}^{S_{k}''}$. 
 From $\gamma S_{k}''=\beta S_{k}$, we have $\tilde{f}_{1}^{S_{k}''}(\gamma)=\tilde{f}_{1}^{S_{k}}(\beta)=\lambda_{\beta}^{S_{k}}$. 
 Furthermore, from $\alpha\gamma S_{k}''=\alpha\beta S_{k}$, 
 we have $\tilde{f}_{1}^{S_{k}''}(\alpha\gamma)=\tilde{f}_{1}^{S_{k}}(\alpha\beta)=\lambda_{\alpha}^{S_{k}}+\lambda_{\beta}^{S_{k}}$. 
 Since $\tilde{f}_{1}^{S_{k}''}$ is modular on $N \setminus S_{k}''$, 
 we obtain 
 \begin{equation}
  \tilde{f}_{1}^{S_{k}''}(\alpha)
   = \tilde{f}_{1}^{S_{k}''}(\alpha\gamma) + \tilde{f}_{1}^{S_{k}''}(\emptyset) - \tilde{f}_{1}^{S_{k}''}(\gamma)
   =\lambda_{\alpha}^{S_{k}}. 
   \label{eq:f1Skdd_modular}
 \end{equation}
 From (\ref{eq:f1Skd_modular}), (\ref{eq:f1Skdd_modular}) and  $\beta S_{k}'=\alpha S_{k}''$ , 
 we have $ \lambda_{\beta}^{S_{k}} = \tilde{f}_{1}^{S_{k}'}(\beta) = \tilde{f}_{1}^{S_{k}''}(\alpha)=\lambda_{\alpha}^{S_{k}}$. \\
\end{proof}

%
For stating other results on skeletal  functions, 
we let $N=\{1,\dots,n\}$, $n=|N|$, and consider $f$ as a 
function of $x_1, \dots, x_n$, where $x_i \in \{0,1\}$, $i=1,\dots,n$.
$S\subseteq N$ can be identified with $(x_1, \dots, x_n)\in \{0,1\}^n$, where
$x_i=1$ if and only if $i\in S$.  Let $\Delta_i f$ denote
the function depending on $x_j$, $j\neq i$, defined as
\[
\Delta_i f= f(x_1, \dots, x_{i-1},1, x_{i+1},\dots, x_n) - 
f(x_1, \dots, x_{i-1},0, x_{i+1},\dots, x_n).
\]
Then it is easily shown that $f$ is supermodular if and only if 
\[
\Delta_i \Delta_j f \ge 0, \quad \forall i\neq j.
\]
Note that $\Delta_i \Delta_j f$ is a function of $x_k$, $k\neq i,j$.
Furthermore $f$ does not depend on $x_i$ if and only if $\Delta_i f\equiv 0$.
Also since a  standardized supermodular function $f$ is non-decreasing,
$\Delta_i f \ge 0$.

From now on we interchangeably write $f(S)$ or $f(x_1, \dots, x_n)$, i.e.,
the argument of $f$ may be a set or a 0-1 vector, depending on the context.
We now prove the following lemma.
\begin{lemma}
\label{lem:marginal-function}
Let $A\subseteq N$ and let $g: {\cal P}(A)\rightarrow \RR$.  
Define $f:{\cal P}(N)\rightarrow \RR$ by $f(S)=g(S\cap A)$.
Then $f$ is supermodular if and only if $g$ is supermodular.  Furthermore
$f$ is skeletal if and only if $g$ is skeletal.
\end{lemma}
\begin{proof}
We omit the proof of the first statement, since it is trivial.
Also it is easily shown that if $g$ is not skeletal, then $f$ is not skeletal.
It remains to show that if $g$ is skeletal, then $f$ is skeletal.
Write $f=f_1 + f_2$, where $f, f_1, f_2$ are standardized.
Then  $0=\Delta_j f$ for $j\not\in A$ implies $0=\Delta_j f_1 =\Delta_j f_2$, $j\not\in A$.
Hence $f_1, f_2$ do not depend on $x_j, j\not\in A$ and $f_1, f_2$ can be considered as
supermodular functions on ${\cal P}(A)$.  Then since $g$ is skeletal, $f_1, f_2$ have
to be proportional to $g=f$.  
\end{proof}

For $A\subseteq N$ let
\begin{equation}
\label{eq:indicator-A}
f_A(S)=1_{A\subseteq \bullet}(S) = 1_{\{A\subseteq S\}} 
=\begin{cases} 1 & \text{ if }  A \subseteq S\\
                 0 & \text{ otherwise.}
\end{cases}
\end{equation}
Using above lemmas we prove the following proposition.
\begin{proposition}
\label{prop:indicator-skeletal}
$f_A = 1_{A\subseteq\bullet}$ is skeletal.
\end{proposition}
\begin{proof}
On ${\cal P}(A)$, $g(S)=g_A(S)=\max(|S|-|A|+1,0)$ is skeletal by 
Lemma \ref{lem:max-k}. Then its extension to ${\cal P}(N)$
is skeletal by Lemma
\ref{lem:marginal-function}.
\end{proof}

We now consider more complicated extensions of
supermodular functions. For $f: \{0,1\}^n \rightarrow \RR$
define
\begin{align*}
f_0(x_1, \dots, x_{n-1})&=f(x_1, \dots, x_{n-1},0), \\
f_1(x_1, \dots, x_{n-1})&=f(x_1, \dots, x_{n-1},1).
\end{align*}
The following lemma is easy to prove.
\begin{lemma}
\label{lem:f0f1}
Let $f$ be a standardized supermodular function, such
that $f_0\equiv 0$.  Then $f$ is skeletal if and only if $f_1$ is skeletal.
\end{lemma}

The following result is somewhat more difficult to prove. 
In the lemma and its proof we denote $N'=\{1,\dots,n-1\}$.
\begin{lemma}
\label{lem:f0f1-2}
Consider $f$ such that $f_0$ is a skeletal function 
satisfying $\Delta_i f_0 (N'\setminus i)= 1$ for all $i\in N'$ and
$f_1(S)=|S|$, $S\subseteq N'$. This 
$f$ is skeletal.
\end{lemma}
\begin{proof}
We first show that $f$ is supermodular.
It suffices to show that $\Delta_i \Delta_j f \ge 0$ for all $i,j \in N, i \not= j$.
For $i,j \not= n$ and $S \subseteq N' \setminus ij$, $\Delta_i \Delta_j f(S) = \Delta_i \Delta_j f_0(S) \ge 0$
because of the supermodularity of $f_0$. Also, for $i,j \not= n$ and $n \in S \subseteq N \setminus ij$,
$\Delta_i \Delta_j f(S\setminus n) = \Delta_i \Delta_j f_1(S) \ge 0$. Thus we need to show the supermodularity when $j = n$ and
$S \subseteq N'$. 
From the supermodularity of $f_{0}$, for $i \in N' \setminus S$ we have 
\[
f_0(Si) - f_0(S) \le 
 f_0(N') - f_0(N' \setminus i) = 1 = f_1(Si) - f_1(S),
\]
which means that 
\[ \Delta_i \Delta_j f(S) = \Delta_i (f_1(S) - f_0(S)) = \Delta_i f_1(S) - \Delta_i f_0(S) \ge 0. \]
Therefore $f$ is supermodular.

Now write $f = g+h$ where $g,h$ are standardized.
Define $g_0, g_1, h_1, h_0$ as above.
Considering the case $x_n=0$, since $f_0$ is skeletal, there exist
$c_g, c_h\ge 0$, $c_g + c_h=1$, such that
\[
g_0 = c_g f_0, \quad h_0 = c_h f_0.
\]
Also from the assumption and supermodularity of $g_{1}$ and $h_{1}$, for any $i\in N'$ we have
\begin{align}
g_1(N') - g_1(N'\setminus i) 
&\ge g_0(N') - g_0(N'\setminus i)= c_g \Delta_i f_0(N'\setminus i)=c_g,
\label{eq:g01}
\\
h_1(N') - h_1(N'\setminus i) 
&\ge h_0(N') - h_0(N'\setminus i)= c_h \Delta_i f_0(N'\setminus i)=c_h.
\label{eq:h01}
\end{align}
If at least one of the inequalities is strict, we have
\begin{align*}
1 &= (n-1)-(n-2)= f_1(N')-f_1(N'\setminus i) \\
&= ( g_1(N') - g_1(N'\setminus i)) + 
 ( h_1(N') - h_1(N'\setminus i))\\
&>  c_g \Delta_i f_0(N'\setminus i) + 
  c_h \Delta_i f_0(N'\setminus i) = 1, 
\end{align*}
which is a contradiction.   Therefore \eqref{eq:g01} and \eqref{eq:h01}
are equalities.

Now consider the case that $x_n=1$.  $f_1(S)=|S|$ is modular. Hence
both $g_1$ and $h_1$ are modular as well.  
From~(\ref{eq:g01}) and the modularity of $g_{1}$, for any $i \in N'$ we have 
\[
c_g = g_1(N')-g_1(N'\setminus i) = 
g_1(i)-g_1(\emptyset)=g_1(i). 
\]
Furthermore from $f_{1}(\emptyset)=0$ and the modularity of $g_{1}$, we have $g_{1}(\emptyset)=0$. 
Hence by (\ref{eq:modular}), we obtain
\[
g_1(S)
 = \sum_{i \in S}g_{1}(i)
 = c_g |S| = c_g f_1(S), \ S\subseteq N'.
\]
Therefore $g=c_g f$.  Similarly $h=c_h f$.  This proves that $f$ is skeletal.
\end{proof}

For our next result, we write ${\bm x}_{n-1}=(x_1, \dots, x_{n-1})$.
Let $f' : \{0,1\}^n \rightarrow \RR$ be supermodular.  Define
$f: \{0,1\}^{n+1} \rightarrow \RR$ by
\[
f({\bm x}_{n-1}, x_n, x_{n+1}) = 
\begin{cases}
f'({\bm x}_{n-1},1) & x_n=x_{n+1}=1, \\
f'({\bm x}_{n-1},0) & \text{otherwise}.
\end{cases}
\]
Then the following result holds.
\begin{proposition}
\label{prop:fig11-extension}
If $f'$ is skeletal, then $f$ is skeletal.
\end{proposition}
\begin{proof}
As in the last proposition we write $N'=\{1,\dots, n-1\}$.  We also write
$N''=\{1,\dots,n+1\}$.
We first show that $f$ is supermodular, i.e.\ 
$\Delta_i \Delta_j f \ge 0$ for all $i,j\in N''$, $i\neq j$.
By assumption,  $\Delta_i \Delta_j f \ge 0$ for $i,j\in N'$.
Next consider the case $i \in N'$ and $j=n+1$. Then by definition
\begin{align}
\label{eq:n11}
\Delta_{n+1} f({\bm x}_{n-1},1)&= \Delta_n f'({\bm x}_{n-1}),\\
\label{eq:n10}
\Delta_{n+1} f({\bm x}_{n-1},0)&= f'({\bm x}_{n-1},0)-f'({\bm x}_{n-1},0)=0.
\end{align}
Hence 
\[
\Delta_i \Delta_{n+1} f({\bm x}_{n-1,\hat{i}},1)=\Delta_i \Delta_n f'({\bm x}_{n-1,\hat{i}})\ge 0, \qquad
\Delta_i \Delta_{n+1} f({\bm x}_{n-1,\hat{i}},0)=0, 
\]
where ${\bm x}_{n-1,\hat{i}}={\bm x}_{N'\setminus i}$. 
The case of $i\in N', j=n$ is similar.
For $i=n, j=n+1$, by taking the difference of 
\eqref{eq:n10} and \eqref{eq:n11}  we have
\[
\Delta_n \Delta_{n+1} f({\bm x}_{n-1})=\Delta_n f'({\bm x}_{n-1})\ge 0.
\]
Therefore $f$ is supermodular.

Now let $f= g + h$.  Repeating $x_n$ define
\[
\tilde f(x_1,\dots, x_n)=f(x_1, \dots, x_n, x_n).
\]
Then $\tilde f=f'$.
By assumption $\tilde f$ is skeletal.  Hence defining $\tilde g, \tilde h$ similarly,
we have
\[
\tilde g=c_g \tilde f, \ \tilde h = c_h \tilde f, \quad c_g, c_h\ge 0, \ 
 c_g + c_h=1.
\]
This implies that $g(x)=c_g f(x)$ and $h(x)=c_h f(x)$ if $x_n=x_{n+1}$.
Therefore we only need to consider the cases $(x_n, x_{n+1})=(1,0)$ or 
$(x_n, x_{n+1})=(0,1)$.  By symmetry we consider only  the former case.
Then
\[
c_g \tilde f ({\bm x}_{n-1},0)=\tilde g({\bm x}_{n-1},0)=g({\bm x}_{n-1},0,0)
\le g({\bm x}_{n-1},1,0).
\]
The last inequality holds by monotonicity of $g$.  Similar relation holds for $h$.
If at least one of the inequalities for $g$ and for $h$ is strict, then
\begin{align*}
\tilde f({\bm x}_{n-1},0)&= c_g \tilde f({\bm x}_{n-1},0)
 +  c_h \tilde f({\bm x}_{n-1},0) \\
&< g({\bm x}_{n-1},1,0)
 +  h({\bm x}_{n-1},1,0) =f({\bm x}_{n-1},1,0) \\
&=f'({\bm x}_{n-1},0)=\tilde f({\bm x}_{n-1},0),
\end{align*}
which is a contradiction. Hence
\[
c_g f({\bm x}_{n-1},1,0)=c_g f'({\bm x}_{n-1},0)=c_g \tilde f({\bm x}_{n-1},0)=g({\bm x}_{n-1},1,0)
\]
and $g=c_g f$. Similarly $h= c_h f$.  Hence $f$ is skeletal.
\end{proof}

Our final result of this section concerns a product of
skeletal functions for disjoint sets.  These functions
played an important role in Kashimura et al\cite{kashimura-etal-2by2}.

\begin{proposition}
\label{prop:product-skeletal}
Let $A B=N$, $A\cap B=\emptyset$ and
consider $f$ of the form
\begin{equation}
\label{eq:multplicative-skeletal}
f(S)=g(A\cap S) h(B\cap S)
\end{equation}
where $g: {\cal P}(A)\rightarrow \RR$ and
$h: {\cal P}(B)\rightarrow \RR$.  Suppose that
$f, g, h$ are standardized supermodular functions.
Then $f$ is skeletal if and only if $g$ and $h$ are skeletal.
\end{proposition}

Note that if $g$ and $h$ are standardized supermodular functions on 
${\cal P}(A)$ and ${\cal P}(B)$, respectively, then $f(S)$ 
in \eqref{eq:multplicative-skeletal} is a standardized supermodular function
on ${\cal P}(N)$.

\begin{proof}
It is easily shown that if $g$ or $h$ are not skeletal, then $f$ is not
skeletal.  Hence if $f$ is skeletal, then both $g$ and $h$ are skeletal.

We now show the converse.  Suppose that $g$ and $h$ are skeletal.
Write $f= f_1 + f_2$, where 
$f_1, f_2$ are standardized.
Consider any $A'\subseteq A$ such that  $g(A')> 0$. 
For $S\subseteq B$, consider
\[
f(A'S)=f_1(A'S) + f_2(A'S) = g(A') h(S), \quad S\subseteq B.
\]
Since $f_1(A'S)$, $f_2(A'S)$ as functions of $S$ are supermodular and since
$h$ is skeletal, there exist $c_1(A'), c_2(A')\ge 0$, such that
\[
f_1(A'S)=c_1(A') h(S),  \ f_2(A'S)=c_2(A')h(S).
\]
Hence 
\[
f(A'S)= c_1(A')h(S)+ c_2(A') h(S) = g(A') h(S), \quad S\subseteq B.
\]
Note that this holds also for $A'$ such that $g(A')=0$ by defining
$0=c_1(A')=c_2(A')$. 
Now fixing any $S\subseteq B$ such that $h(S)>0$ 
and considering
$f(A'S)/h(S)=g(A')=c_1(A')+c_2(A')$ as a function of $A'$,
we see that $c_1$ and $c_2$ have to 
be proportional to $g$ because $g$ is skeletal. 
Then $f_1$ and $f_2$ are proportional to $f=gh$.
\end{proof}

In Appendix A, we present a generalization of Proposition 4.3 for
arbitrary cones.

\section{Structure of faces of semi-elementary imsets}
\label{sec:semi-elementary-face}
\newcommand{\cEABC}{\cE_{\tri{A}{B}{C}}}
\newcommand{\tabG}{{\tri{a}{b}{\Gamma}}}

\newcommand{\FABC}{F_{\tri{A}{B}{C}}}
\newcommand{\MABC}{{\cal M}_{\tri{A}{B}{C}}}
\newcommand{\subind}[1]{1_{#1 \subseteq \bullet}}
\newcommand{\supind}[1]{1_{#1 \supseteq \bullet}}

In this section we study a face of the imset cone corresponding to a
semi-elementary imset $u_{\tri{A}{B}{C}}$.  By the face corresponding
to $u_{\tri{A}{B}{C}}$, we mean the unique face of the imset cone, such
that $u_{\tri{A}{B}{C}}$ is in the relative interior of the face.
We denote this face by $\FABC$. 
In Kashimura et al.\cite{kashimura-etal-2by2} we have shown 
some remarkable facts on $\FABC$.
Here we establish more basic facts on $\FABC$.

If we repeatedly apply the decomposition in \eqref{eq:sg-imset},
any semi-elementary imset $u_{\tri{A}{B}{C}}$ can be written
as a sum of elementary imsets.  Those elementary imsets are from the
following set of elementary imsets:
\[
\cEABC = \{u_\tabG \mid a\in A, b\in B, 
C\subseteq \Gamma\subseteq ABC\}.
\]
In this section we establish the following basic facts on
$\FABC$:
i) $\cEABC$ is the set of extreme rays of $\FABC$,
ii) $\dim(\FABC)=(2^{|A|}-1)(2^{|B|}-1)$.
For proving these facts we give linearly independent set
of supermodular functions orthogonal to $\FABC$.

Let $D=N\setminus(ABC)$, i.e.\ $N=ABCD$.
For $E\subseteq N$ define $\subind{E}$ by
\eqref{eq:indicator-A}, which is skeletal by Proposition 
\ref{prop:indicator-skeletal}.  Similarly define
\[
\supind{E}(S) = 
\begin{cases} 1 & \text{ if }  E \supseteq S\\
                 0 & \text{ otherwise},
\end{cases}
\]
which is skeletal by Lemma \ref{lem:complementation}, since
$\supind{E} = \subind{E^C}^C$.
Now consider the following supermodular functions.
\begin{align}
 & \subind{A_1C}\quad (A_1\subseteq A),
 \\
 & \subind{B_1C}\quad (\emptyset \neq B_1\subseteq B),
 \\
 & \supind{EC_1}\quad (E\subseteq AB,\ C_1\subset C),
 \label{eq:ortho}
 \\
 & \subind{ED_1}\quad (E\subseteq ABC,\ \emptyset\neq D_1\subseteq D).
\end{align}
Let $\MABC$ denote the set of these supermodular functions.
The cardinality of $\MABC$ is given by
\begin{align}
 |\MABC|&=2^{|A|}+(2^{|B|}-1)+2^{|AB|}(2^{|C|}-1)+2^{|ABC|}(2^{|D|}-1)
 \nonumber \\
 &= 2^{|N|} - (2^{|A|}-1)(2^{|B|}-1).
\label{eq:cardinality-MABC}
\end{align}

First we check that the linear independence of the above
supermodular functions.

\begin{lemma} The elements of $\MABC$ are linearly independent.
\label{lem:MABC-independent}
\end{lemma}
\begin{proof}
We need to show that the coefficients 
$\lambda_{A_1},\mu_{B_1},\nu_{EC_1},\xi_{ED_1}\in\mathbb{R}$
in 
\begin{align*}
 &\sum_{A_1\subseteq A}\lambda_{A_1}\subind{A_1C}(S) 
 +\sum_{\emptyset\neq B_1\subseteq B}\mu_{B_1} \subind{B_1C}(S) 
 +\sum_{E\subseteq AB,C_1\subset C}\nu_{EC_1} \supind{EC_1}(S) 
 \\
 &\quad +\sum_{E\subseteq ABC,\emptyset\neq D_1\subseteq D}\xi_{ED_1} \subind{ED_1}(S) 
 = 0\quad  (\forall S)
\end{align*}
are zeros.
Let $S=A_2C$ ($A_2\subseteq A$).  Then 
$\sum_{A_1\subseteq A_2}\lambda_{A_1}=0$ and by induction we have
$\lambda_{A_1}=0$.
Similarly by letting $S=B_2C$ ($\emptyset\neq B_2\subseteq B$) we have
$\sum_{\emptyset\neq B_1\subseteq B_2}\mu_{B_1}=0$ and by induction 
$\mu_{B_1}=0$.

Next let $S=E_2C_2$ ($E_2\subseteq AB,C_2\subset C$). Then
\[
\sum_{E_2\subseteq E\subseteq AB,C_2\subseteq C_1\subset C}\nu_{EC_1}=0
\]
and the double induction on $E$ and $C_1$ yields
$\nu_{E C_1}=0$.  Finally 
letting $S=E_2D_2$ ($E_2\subseteq ABC,\emptyset\neq D_2\subseteq D$)
we have 
$\sum_{E\subseteq E_2,\emptyset\neq D_1\subseteq D_2}\xi_{ED_1}=0$
and the double induction on $E$ and $D_1$ yields $\xi_{ED_1}=0$.
\end{proof}

Now we prove the following theorem.
\begin{theorem}
\label{thm:semi-elementary-orthogonal}
An elementary imset $u$ belongs to  $\cEABC$ if and only if
\begin{equation}
\label{eq:ABC-orthogonality}
\langle u, f \rangle = 0, \quad \forall f\in \MABC.
\end{equation}
\end{theorem}

\begin{proof}
We first show the orthogonality 
$\langle u, f \rangle = 0$ for all $u\in \cEABC$ and for all $f\in \MABC$.
Consider $f(S)=\subind{A_1C}$ 
($A_1\subseteq A$). 
Let 
\[
 u_{\tri{a}{b}{C'}}=\delta_{abC'}+\delta_{C'} - \delta_{aC'}-\delta_{bC'}
\in \cEABC.
 \]
Then 
\begin{align*}
\langle  u_{\tri{a}{b}{C'}}, f\rangle
 &=f(abC')+f(C')-f(aC')-f(bC')
 \\
 &= 1_{\{A_1C  \subseteq abC' \}} + 1_{\{ A_1C \subseteq C' \}}
 -1_{\{A_1C  \subseteq aC'\}} - 1_{\{ A_1C \subseteq bC' \}}
 \\
 &= 1_{\{ A_1C \subseteq  aC' \}} + 1_{\{ A_1C \subseteq  C'\}}
 -1_{\{ A_1C \subseteq aC' \}} - 1_{\{ A_1C \subseteq C'\}}
 \\
 &= 0.
\end{align*}
Similarly $f(S)=\subind{B_1C}$ 
is orthogonal to 
$u_{\tri{a}{b}{C'}}\in \cEABC$.

Next consider $f(S)=\supind{EC_1}$ 
($E\subseteq AB,C_1\subset C$).
Since $C\not\subseteq C_1$
\begin{align*}
\langle  u_{\tri{a}{b}{C'}}, f\rangle 
 &= f(abC')+f(C')-f(aC')-f(bC')
 \\
 &= 1_{\{ EC_1 \supseteq abC' \}} + 1_{\{ EC_1  \supseteq C'\}}
 - 1_{\{  EC_1\supseteq aC' \}} - 1_{\{EC_1  \supseteq bC' \}}
 \\
 &= 0+0-0-0=0.
\end{align*}
Finally consider $f(S)=\subind{ED_1}$ 
($E\subseteq ABC,
\emptyset\neq D_1\subseteq D$). Then
\begin{align*}
\langle  u_{\tri{a}{b}{C'}}, f\rangle
 &=f(abC')+f(C')-f(aC')-f(bC')
 \\
 &= 1_{\{ ED_1 \subseteq abC' \}} + 1_{\{ ED_1 \subseteq C' \}}
 - 1_{\{ED_1 \subseteq aC'\}} - 1_{\{ ED_1 \subseteq bC' \}}
 \\
 &= 0 + 0 - 0 - 0=0.
\end{align*}
We have shown that 
$\langle u, f \rangle = 0$, $\forall u\in \cEABC$, $\forall f\in \MABC$.

We now show that for any elementary imset $u=u_{\tri{a'}{b'}{C'}}\not\in \cEABC$ there exists
$f\in \MABC$ such that $\langle  u, f\rangle=1$.
For $f=\subind{F}$ we here confirm which elementary imset
$u_{\tri{\alpha}{\beta}{\Gamma}}$ satisfies 
$\langle u_{\tri{\alpha}{\beta}{\Gamma}}, f \rangle=1$.
Write
\[
\Gamma_1 = \Gamma\cap F, \quad \Gamma_2 = \Gamma\cap F^C.
\]
Then we can ignore $\Gamma_2$ and replace $\Gamma$ by $\Gamma_1$ in
\begin{equation}
\label{eq:1111}
1_{F\subseteq \alpha\beta\Gamma} + 1_{F\subseteq \Gamma}
 - 1_{F\subseteq \alpha\Gamma} - 1_{F\subseteq \beta\Gamma}.
\end{equation}
Hence we only consider $\Gamma_1$.  It is easy to see that
\eqref{eq:1111} is equal to 1 if and only if the first term
of \eqref{eq:1111} is 1 and other three terms are zeros. In this case
\[
|\alpha\beta\Gamma_1| \ge |F|, \quad |\Gamma_1| \le |F|-2
\]
implies  $F=\alpha\beta\Gamma_1$.  Putting $\Gamma_2$ back we have the following equivalence
\begin{align}
\langle u_{\tri{\alpha}{\beta}{\Gamma}}, \subind{F}\rangle = 1
&\ \Leftrightarrow \ \alpha, \beta\in F, \ (\alpha\beta\Gamma)\cap F=F 
\nonumber \\
&\ \Leftrightarrow \ \alpha, \beta\in F,\ \alpha\beta\Gamma \supseteq F.
\label{eq:non-orthogonality-with-subind}
\end{align}
Similarly for $\supind{F}$  we can show
\begin{equation}
\langle u_{\tri{\alpha}{\beta}{\Gamma}}, \supind{F}\rangle = 1
\ \Leftrightarrow \ \alpha, \beta\not\in F,\ \Gamma \subseteq F.
\label{eq:non-orthogonality-with-supind}
\end{equation}

Now arbitrarily fix $\alpha,\beta,\Gamma$. Consider the case
$\alpha\beta\Gamma\not\subseteq ABC$, i.e.\ $(\alpha\beta\Gamma)\cap D \neq \emptyset$. 
Let
\[
D_1 = (\alpha\beta\Gamma)\cap D, \ E=(\alpha\beta\Gamma) \cap ABC.
\]
Then $ED_1 = \alpha\beta\Gamma$ and 
$\alpha,\beta\in F$, $\alpha\beta \Gamma \supseteq F$, where 
$F=ED_1$.
Hence (\ref{eq:non-orthogonality-with-subind}) implies 
$
\langle u_{\tri{\alpha}{\beta}{\Gamma}}, \subind{ED_1} \rangle=1 .
$
It remains to consider the cases where $\alpha\beta\Gamma\subseteq ABC$.
Hence assume $C\not\subseteq \Gamma\subseteq ABC$.
Let $C_1 =C\cap \Gamma$.  Then $C_1 \subset C$. 
Let $E=\Gamma \cap (AB)$.  Then $\Gamma=EC_1$.  For $F=\Gamma$,
we have
$
\alpha,\beta\not\in F, \Gamma=F.
$
Hence \eqref{eq:non-orthogonality-with-supind} implies
$
\langle u_{\tri{\alpha}{\beta}{\Gamma}}, \supind{EC_1} \rangle=1 .
$
It remains to consider the case
$
\Gamma\supseteq C, \ \alpha\beta\Gamma\subseteq ABC.
$
Note that we want to eliminate the case that both $\alpha$ and $\beta$
belong to $A$ or belong to $B$.  By symmetry we can only consider the
former case.
Hence assume 
$\alpha,\beta\in A, \ \Gamma\supseteq C$. 
In this case we can set $A_1 = \alpha\beta (A\cap \Gamma)$.
Then (\ref{eq:non-orthogonality-with-subind}) implies 
$
\langle u_{\tri{\alpha}{\beta}{\Gamma}}, \subind{A_1 C} \rangle=1 .
$
\end{proof}

\begin{remark} In the above proof we have shown that
\label{rem:0or1}
\[
\langle u, \subind{A}\rangle, 
\langle u, \supind{A}\rangle \in \{0,1\}
\]
for every elementary imset $u$ and every $A\subseteq N$.
\end{remark}

\begin{remark}
\label{rem:ABC-linear-orthogonality}
Let ${\rm Lin}(\MABC)$ and ${\rm Lin}(\cEABC)$ denote
linear subspaces spanned by $\MABC$ and $\cEABC$, respectively.
Linearly extending  \eqref{eq:ABC-orthogonality} we see that
${\rm Lin}(\MABC)$ is the orthogonal complement of 
${\rm Lin}(\cEABC)$.
\end{remark}

Two facts stated at the beginning of this section are simple
consequences of Theorem  \ref{thm:semi-elementary-orthogonal}.
\begin{corollary}
$\cEABC$ is the set of extreme rays of $\FABC$ and 
\[
\dim(\FABC)=(2^{|A|}-1)(2^{|B|}-1).
\]
\end{corollary}
The first statement is obvious and the second statement follows from
Remark \ref{rem:ABC-linear-orthogonality} and \eqref{eq:cardinality-MABC}.
The fact that  $\dim(\FABC)$ does not depend on $C$ can also be seen from
the one-to-one linear correspondence between ${\rm Lin}(\cEABC)$ and 
${\rm Lin}(\cE_{\tri{A}{B}{\emptyset}})$ given by
\[
\cEABC\ni u_{\tri{a}{b}{\Gamma}}\  \leftrightarrow \ u_{\tri{a}{b}{\Gamma\setminus C}}\in
\cE_{\tri{A}{B}{\emptyset}}.
\]

Note that Theorem \ref{thm:semi-elementary-orthogonal}
shows that elements of $\MABC$ are extreme rays of the
face $\FABC^\perp$ 
of the supermodular cone $K(N)$ dual to $\FABC$.
However it is hard to believe that $\FABC^\perp$ is simplicial
and hence $\FABC^\perp$ should have other extreme rays than
elements of $\MABC$.  Characterization of extreme rays of 
$\FABC^\perp$ is an interesting question left to our future research.

\begin{remark}
\label{rem:2by2-comment}
In our previous manuscript Kashimura et al\cite{kashimura-etal-2by2}.
we studied how the semi-elementary imset $u_{\tri{A}{B}{C}}$ is expressed as a
non-negative integer combination of elements of $\cEABC$.  By 
Theorem \ref{thm:semi-elementary-orthogonal}, if $u_{\tri{A}{B}{C}}$ is expressed
as a non-negative integer combination of all elementary imsets, then the coefficients of 
$u\not\in \cEABC$ have to be zero.  Therefore  Theorem \ref{thm:semi-elementary-orthogonal}
justifies our restriction to imsets from $\cEABC$ for expressing $u_{\tri{A}{B}{C}}$.
\end{remark}

\section{Small  relations among imsets}
\label{sec:small-linear-relations}
Consider the configuration $\cfg_N$ as in Table
\ref{tab:configuration}.  The kernel of $\cfg_N$ is denoted by $\ker \cfg_N$.
Furthermore let
\[
\ker_{\ZZ} \cfg_N = (\ker \cfg_N) \cap \ZZ^{|\cE(N)|}
\]
denote the integer kernel of $\cfg_N$.
Let $z=(z_u)_{u\in \cE(N)} 
\in  \RR^{|\cE(N)|}$ 
be an element of the kernel of $\cfg_N$.
Considering $f$ in \eqref{eq:homogeneous-configuration} as a row vector
we have
\begin{equation}
\label{eq:relation-sum0}
0 = f \cfg_N z = (1,1,\dots,1) z = \sum_{u\in \cE(N)} z_u.
\end{equation}
Hence any non-zero $z\in \ker \cfg_N$ has both positive and
negative elements. We denote the positive elements of $z$
as $\alpha_1, \dots, \alpha_k$ and negative elements of $z$ as
$-\beta_1, \dots, -\beta_m$.  
Then \eqref{eq:relation-sum0}
can be written as
\begin{equation}
\label{eq:relation}
\alpha_1 u_1 + \dots + \alpha_k u_k = \beta_1 v_1 + \dots +  \beta_m v_m, 
\end{equation}
where $u_1,\dots, u_k, v_1, \dots, v_m$ are distinct elementary imsets.
We call this kind of equality between two non-negative combinations
of elementary imsets a {\it relation} among elementary imsets.
More precisely we call \eqref{eq:relation} a $k$ by $m$ relation.
Since $-z$ also belongs to $\ker \cfg_N$, 
we can take  $k\le m$. 
Also we can assume that $u_1,\dots,u_k$ and $v_1,\dots, v_m$
are ordered according to the graded reverse lexicographic order of elementary imsets. 
We say that the  relation \eqref{eq:relation} contains another relation
$\sum_{i'=1}^{k'} \alpha'_{i'} u'_{i'} = \sum_{j'=1}^{m'} \beta'_{j'} v'_{j'}$ 
if every $u'_{i'}$ appears on the left-hand side of \eqref{eq:relation}
or every $v'_{j'}$ appears on the right-hand side of \eqref{eq:relation}.

The most basic relation is the following 2 by 2 relation, 
which comes from the semi-graphoid
axiom.
\begin{equation}
\label{eq:2:2-relation}
u_{\tri{a}{b_1}{C}} + u_{\tri{a}{b_2}{b_1C}}
=u_{\tri{a}{b_2}{C}} + u_{\tri{a}{b_1}{b_2 C}}.
\end{equation}
The corresponding element of $\ker \cfg_N$ is written as
\begin{equation}
\label{eq:2-2-relation-kernel}
 \delta_{\tri{a}{b_1}{C}} + \delta_{\tri{a}{b_2}{b_1C}}
 - \delta_{\tri{a}{b_2}{C}} - \delta_{\tri{a}{b_1}{b_2C}},
\end{equation}
where the unit vector $\delta_{\tri{a}{b}{C}}\in \ZZ^{|\cE(N)|}$
is defined as $\delta_{\tri{a}{b}{C}}(v)=1$ if $v=u_{\tri{a}{b}{C}}$
and $0$ otherwise.  In connection to Markov bases discussed in Section \ref{sec:markov-basis}
we call an element of $\ker \cfg_N$  a {\it move}.
Let $\cB$ denote the collection of  moves in \eqref{eq:2-2-relation-kernel}
over all $a,b_1,b_2,C$. We first show that $\cB$ generates the integer kernel of $\cfg_N$,
i.e.\ $\cB$ contains a lattice basis for $\cfg_N$.

\begin{theorem}
\label{thm:kernel-lattice-basis}
$\ker_\ZZ \cfg_N$ as a submodule of $\ZZ^{|\cE(N)|}$ is generated by $\cB$.  
\end{theorem}

The following proof also shows that $\cB$ spans $\ker \cfg_N$ as a linear subspace of 
$\RR^{|\cE(N)|}$.

\begin{proof}
 Let $N=\{1,\ldots,n\}$.
 Let $\ZZ\cB$ be the submodule of $\ZZ^{|\cE(N)|}$ generated by $\cB$.
 Since the inclusion $\ker_\ZZ \cfg_N \supseteq \ZZ\cB$ holds,
 we prove $\ker_\ZZ \cfg_N\subseteq \ZZ\cB$.
 Consider the graded reverse lexicographic order of $\cE(N)$ in \eqref{eq:grlex-for-elementary}.
 Recall that we always assume $a<b$ if we write $u_{\tri{a}{b}{C}}$
 and the maximum element of $\cE(N)$ with respect to $\precw$ is $u_{\tri{1}{2}{R}}$, where $R=\{3,\ldots,n\} $.

 For any non-zero $z\in \ZZ^{|\cE(N)|}$
 we define the degree $\deg(z)$ by the {\em minimum} elementary imset $u_{\tri{a}{b}{C}}$
 such that $z(u_{\tri{a}{b}{C}})\neq 0$.
 Now fix any non-zero $z\in\ker_\ZZ \cfg_N$.
 We will find some $v\in\ZZ\cB$
 such that $\deg(z-v)\succw\deg(z)$ or $z=v$.
 Then the proposition follows by induction.
 We first show that $\deg(z)$ is not the maximum $u_{\tri{1}{2}{R}}$.
 Indeed, if $\deg(z)=u_{\tri{1}{2}{R}}$,
 then $z$ has to be $z(u_{\tri{1}{2}{R}})\delta_{\tri{1}{2}{R}}$ and therefore
 \begin{align*}
  0 = \cfg_N\cdot z
  = z(u_{\tri{1}{2}{R}})u_{\tri{1}{2}{R}},
 \end{align*}
 which contradicts to $z(u_{\tri{1}{2}{R}})\neq 0$.
 Hence $\deg(z)\precw u_{\tri{1}{2}{R}}$.
 Denote $\deg(z)=u_{\tri{a}{b}{C}}$
 and let $\cE_C=\{u_{\tri{a'}{b'}{C}}\mid a',b'\in N\setminus C\}$.
 Denote the maximum of $\cE_C$ by $u_{\tri{\alpha}{\beta}{C}}$.
 We show that $u_{\tri{a}{b}{C}}\precw u_{\tri{\alpha}{\beta}{C}}$.
 Indeed, if $u_{\tri{a}{b}{C}}=u_{\tri{\alpha}{\beta}{C}}$,
 then $z$ can be written as
 \begin{align*}
  z &= z(u_{\tri{\alpha}{\beta}{C}})\delta_{\tri{\alpha}{\beta}{C}}
 + \sum_{C'>C}\sum_{a',b'\in N\setminus C'}z(u_{\tri{a'}{b'}{C'}}) \delta_{\tri{a'}{b'}{C'}}.
 \end{align*}
 Then
 \begin{align*}
  0 = \cfg_N\cdot z
  &= z(u_{\tri{\alpha}{\beta}{C}})u_{\tri{\alpha}{\beta}{C}}
  + \sum_{C'>C}\sum_{a',b'\in N\setminus C'}z(u_{\tri{a'}{b'}{C'}})u_{\tri{a'}{b'}{C'}}
  \\
  &= z(u_{\tri{\alpha}{\beta}{C}})\delta_C + (\mathrm{terms\ of\ }\{\delta_{C'}\}_{C'>C}),
 \end{align*}
 which contradicts to $z(u_{\tri{\alpha}{\beta}{C}})\neq 0$.
 Hence $u_{\tri{a}{b}{C}}\precw u_{\tri{\alpha}{\beta}{C}}$.

 Now we give a vector $v\in\ZZ\cB$ such that $\deg(z-v)\succw\deg(z)$ or $z=v$.
 If $b<\beta$, then consider a vector
 \begin{align*}
  v=z(u_{\tri{a}{b}{C}})
  \left(
  \delta_{\tri{a}{b}{C}} + \delta_{\tri{a}{\beta}{bC}}
  - \delta_{\tri{a}{\beta}{C}} - \delta_{\tri{a}{b}{\beta C}}
  \right)
 \end{align*}
 in $\ZZ\cB$. Note that $\deg(v)=u_{\tri{a}{b}{C}}=\deg(z)$ and the leading terms of $v$ and $z$ are the same.
 Therefore we have $\deg(z-v)\succw\deg(z)$ or $z-v=0$.
 If $b=\beta$, then $a<\alpha$.
 Consider a vector
 \begin{align*}
  v=z(u_{\tri{a}{b}{C}})
  \left(
  \delta_{\tri{a}{b}{C}} + \delta_{\tri{\alpha}{b}{aC}}
  - \delta_{\tri{\alpha}{b}{C}} - \delta_{\tri{a}{b}{\alpha C}}
  \right)
 \end{align*}
 in $\ZZ\cB$.
 Note that $\deg(v)=u_{\tri{a}{b}{C}}=\deg(z)$ and the leading terms of $v$ and $z$ are the same.
 Therefore we have $\deg(z-v)\succw\deg(z)$ or $z-v=0$.
 This completes the proof.
\end{proof}

For the rest of this section we consider $k=2,3$ in 
\eqref{eq:relation}
and show that
only certain types of relations exist.
By \eqref{eq:relation-sum0} we have
\begin{equation}
\label{eq:relation-sum1}
\alpha_1 + \dots + \alpha_k = \beta_1 + \dots + \beta_m
\end{equation}

Our first result on small relations is the following.
\begin{theorem}
\label{thm:2:2}
Consider a relation in $\eqref{eq:relation}$ with $k=2$: 
\begin{equation}
 \alpha_{1} u_{1} + \alpha_{2} u_{2} = \beta_{1} v_{1} + \dots +  \beta_{m} v_{m}, 
 \label{eq:2-by-m-relation}
\end{equation}
where $\alpha_{1}$, $\alpha_{2}$, $\beta_{1}$, \dots, $\beta_{m}$ are integers. 
Then
$m=2$ and $\alpha_1 = \alpha_2 = \beta_1 = \beta_2$ and 
the relation is a positive multiple of \eqref{eq:2:2-relation}.
\end{theorem}

\begin{proof}
Let $u=\alpha_{1}u_{1}+\alpha_{2}u_{2}$. 
From~(\ref{eq:relation-sum1}), we have 
\begin{equation}
 \alpha_1 + \alpha_2 = \beta_1 + \dots + \beta_m
  , \quad 
  \alpha_1 > 0, \alpha_2 > 0, 
  \beta_1 > 0, \dots, \beta_m > 0. 
 \nonumber 
\end{equation}
Because $u_{1}, u_{2}, v_{1}, \dots, v_{m}$ are extreme rays of $\cK^*(N)$, 
there is no relation like $\alpha_{1} u_{1} + \alpha_{2} u_{2} = \beta_{1} v_{1}$. 
Therefore, we have $m \geq 2$. 
Furthermore, because $u_{1}, u_{2}, v_{1}, \dots, v_{m}$ are distinct elementary imset, 
the following properties  hold: 
\begin{equation}
 \quad u_{1} \neq u_{2},\quad u_{1}, u_{2} \neq v_{1}, \dots, v_{m},\quad  v_{i} \neq v_{j}
 \label{eq:proof_theorem3:1.5}
\end{equation}
for $i \neq j$. 
We write the elementary imsets in (\ref{eq:2-by-m-relation}) as 
$u_{1}=\utri{a_{1}}{b_{1}}{C_{1}}, u_{2}=\utri{a_{2}}{b_{2}}{C_{2}}, v_{1}=\utri{a_{1}'}{b_{1}'}{C_{1}'}, \dots, v_{m}=\utri{a_{m}'}{b_{m}'}{C_{m}'}$. 
Remember that the terms are ordered such that $u_{1} < u_{2}$ and $v_{1} < v_{2} < \dots < v_{m}$ according to the graded reverse lexicographic order of elementary imsets. 
Then there exists $l_{1} \in \{1,\dots,m\}$ such that $C_{1} = C_{1}' = \dots =C_{l_{1}}'$
and $C_{l_{1}+1}',\dots, C_{m}' \neq C_{1}$. 
Note that if $C_{1} \neq C_{2}$, then $\alpha_{1} = \beta_{1} + \dots + \beta_{l_{1}}$ and $a_{1}b_{1}C_{1} \neq a_{i}'b_{i}'C_{i}'$ for $i \in \{1,\dots,l_{1}\}$. 
In the same way, there exists $l_{2} \in \{1,\dots,m\}$ such that 
\begin{equation}
 a_{2}b_{2}C_{2} = a_{l_{2}}'b_{l_{2}}'C_{l_{2}}'
  \label{eq:proof_theorem3:2}
\end{equation}
and $C_{2} \neq C_{l_{2}}'$. 


Now we define the weight $w$ for a relation as follows:
\begin{equation}
 w
  = \sum_{ \{S \subseteq N \: \mid \: \abs{S}=\abs{C_{1}}+1 \} }u(S) 
  = \sum_{ \{S \subseteq N \: \mid \: \abs{S}=\abs{C_{1}}+1 \} }(\alpha_{1}u_{1}(S)+\alpha_{2}u_{2}(S)). 
  \nonumber 
\end{equation}
Then the possible values of $w$ are $-2\alpha_{1}$, $-2\alpha_{1}-2\alpha_{2}$ and $-2\alpha_{1} + \alpha_{2}$. 

Firstly, we consider the case of $w = -2\alpha_{1}$. 
Then we have $\abs{C_{2}} \geq \abs{C_{1}} + 2$. 
In this case, $u_{1} = v_{j}$ holds for some $j \in \{1,\dots,m\}$ and this contradicts~(\ref{eq:proof_theorem3:1.5}). 
Therefore there is no relation with $w = -2\alpha_{1}$. 

Next, we consider the case of $w = -2\alpha_{1}-2\alpha_{2}$. 
We prove $C_{1} \neq C_{2}$. 
If $C_{1} = C_{2}$, then by~(\ref{eq:proof_theorem3:2}) there exists some $l \in \{1,\dots,m\}$ such that $v_{l} = u_{2}$. 
This contradicts (\ref{eq:proof_theorem3:1.5}). 
Hence we obtain $C_{1} \neq C_{2}$. 
In the same way, we can prove $a_{1}b_{1}C_{1} \neq a_{2}b_{2}C_{2}$. 
Now, since $\abs{C_{1}}=\abs{C_{2}}$ in this case, we have $C_{2} = C_{l_{1}+1}' = \dots = C_{m}'$. 
However, there exists $l_{2} \in \{1,\dots,m\}$ such that (\ref{eq:proof_theorem3:2}) holds. 
If $l_{2} \in \{l_{1}+1,\dots,m\}$, then this contradicts~(\ref{eq:proof_theorem3:1.5}). 
If $l_{2} \in \{1,\dots,l_{1}\}$, then $C_{1} \subseteq a_{2}b_{2}C_{2}$ and $C_{2} \subseteq a_{1}b_{1}C_{1}$. 
This means that $C_{1}, C_{2} \subseteq a_{1}b_{1}C_{1} \cap a_{2}b_{2}C_{2}$. 
Then we have $C_{1} \cup C_{2} \subseteq a_{1}b_{1}C_{1} \cap a_{2}b_{2}C_{2}$. 
From $C_{1} \neq C_{2}$ and $\abs{C_{1}}=\abs{C_{2}}$, we have $\abs{C_{1} \cup C_{2}} \geq \abs{C_{1}} + 2$. 
Furthermore, from $a_{1}b_{1}C_{1} \neq a_{2}b_{2}C_{2}$, we have $\abs{a_{1}b_{1}C_{1} \cap a_{2}b_{2}C_{2}} \leq \abs{a_{1}b_{1}C_{1}} - 2 = \abs{C_{1}}$. 
However, this means that $\abs{C_{1}} + 2 \leq \abs{C_{1}}$. This is a contradiction. 

Finally, we consider the case of $w = -2\alpha_{1} + \alpha_{2}$. 
Let $h = \abs{C_{1}}$. 
Then we have $\abs{C_{2}}=h+1$, $\abs{\abC{1}}=h+2$, $\abs{\abC{2}}=h+3$ and $\abs{C_{i}'}=h+1$ for $i \in \{l_{1}+1,\dots,m\}$. 
This means that $\abC{2} = \abCd{l_{1}+1} = \dots = \abCd{m}$. 
Now, assume $u(\abC{1}) \geq 1$. 
Since $v_{i}(S) \leq 0$ for $\abs{S}=h+2$ and $i \in \{l_{1}+1,\dots,m\}$, 
we have $\abC{1} = \abCd{j}$ for some $j \in \{1,\dots,l_{1}\}$. 
This contradicts (\ref{eq:proof_theorem3:1.5}). Therefore $ u(\abC{1}) \leq 0$ and we have
\begin{equation}
 0 \geq u(\abC{1}) = \alpha_{1}u_{1}(\abC{1}) + \alpha_{2}u_{2}(\abC{1}) = \alpha_{1} + \alpha_{2}u_{2}(\abC{1}). 
  \nonumber 
\end{equation}
From the above equation, we have $u_{2}(\abC{1}) \leq -\frac{\alpha_{1}}{\alpha_{2}} < 0$. 
This means $u_{2}(\abC{1}) = -1$. Therefore, we have $\alpha_{2} \geq \alpha_{1}$. 
In the similar way, we obtain $u(C_{2}) \leq 0$ and $\alpha_{1} \geq \alpha_{2}$. 
Hence we have $\alpha_{1} = \alpha_{2}$. 
From this, we have $u(C_{2}) = u(\abC{1}) = 0$. 
Then there exist $d_{1}, d_{2} \in \{a_{1},b_{1}\}$ and $d_{3} \in \{a_{2},b_{2}\}$ such that 
\begin{eqnarray}
 C_{2} = C_{1}\cup\{d_{1}\}, \quad
  \abC{1} = C_{2}\cup\{d_{2}\} = C_{1}\cup\{d_{1}, d_{2}\}, 
 \nonumber \\
  \abC{2} = \abC{1}\cup\{d_{3}\} = C_2\cup\{d_{2}, d_{3}\} = C_{1}\cup\{d_{1}, d_{2}, d_{3}\}. 
   \nonumber 
\end{eqnarray}
Then, from $u(C_{1}\cup\{d_{2}\})=u(\abC{2}\setminus\{d_{2}\})=-1$, 
there exist $i' \in \{1,\dots,l_{1}\}$ and $j' \in \{l_{1}+1,\dots,m\}$ such that 
\begin{equation}
 v_{i'}(C_{1}\cup\{d_{2}\})=v_{j'}(\abC{2}\setminus\{d_{2}\})=-1. 
  \label{eq:proof_theorem3:5}
\end{equation}
Let $v_{i'}(C_{1}\cup d_{4})=-1$ for some $d_{4} \in N\setminus\{C_{1}\cup d_{2}\}$. 
Then there exist $j'$ such that $v_{j'}(C_{1}\cup d_{4})=1$. 
Since $v_{j'}(\abC{2})=1$ for this $j'$, we have $C_{1}\cup d_{4} \subseteq \abC{2}=C_{1}\cup\{d_{1}, d_{2}, d_{3}\}$. 
This means $d_{4} \in \{d_{1}, d_{3}\}$. 
Furthermore, since $v_{i'} \neq v_{j'}$ for $i' \neq j'$ in (\ref{eq:proof_theorem3:1.5}), $v_{i'}$ is limited to one of 
$v_{i'} = u_{\tri{d_{2}}{d_{3}}{C_{1}}}$ or $v_{i'} = u_{\tri{d_{3}}{d_{2}}{C_{1}}}$, 
according to $d_{2} < d_{3}$ or $d_{3} < d_{2}$. 
In the same way, $v_{j'}$ is limited to one of 
$u_{\tri{d_{1}}{d_{2}}{d_{3}C_{1}}}$ or $u_{\tri{d_{2}}{d_{1}}{d_{3}C_{1}}}$, 
according to $d_{1} < d_{2}$ or $d_{2} < d_{1}$. 
Therefore we obtain $m=2$. 
Furthermore from \eqref{eq:proof_theorem3:5}, the relation is a positive multiple of  \eqref{eq:2:2-relation}. 
\end{proof}

One consequence of this theorem is the following corollary on 
two-dimensional and three-dimensional faces of $\cK^*(N)$.

\begin{corollary}
Two-dimensional faces of $\cK^*(N)$ are
simplicial cone generated by two arbitrary elementary imsets which
do not appear on the same side of the relation
\eqref{eq:2:2-relation}.  Three-dimensional faces are
either simplicial cones having three elementary imsets as extreme rays
or a cone with four extreme rays appearing in the relation
\eqref{eq:2:2-relation}.
\end{corollary}
\begin{proof}
Two-dimensional cones are always simplicial.  If two elementary imsets
do not span a two-dimensional face, then the two-dimensional cone spanned by
these two imsets cuts a relative interior of a face of a higher dimension.
Then it intersects a cone generated by other extreme rays of this face.  Hence
these two imsets appear on one side of a  2 by  $m$ relation.  
Then by Theorem \ref{thm:2:2} $m=2$ and this face has to be of the form \eqref{eq:2:2-relation}.
If a three-dimensional cone is not simplicial, then any set of
four extreme rays are linearly dependent.  Hence there is a relation
among them, which has to be a two by two relation.  By the same argument it is easily seen 
that there are no more than four extreme rays of a three-dimensional face.
Then by Theorem \ref{thm:2:2} this face has to be
of the form \eqref{eq:2:2-relation}.
\end{proof}

Another important example of relation is \eqref{eq:3:3-relation} with $k=m=3$.

\begin{theorem}
\label{thm:3:3}
Consider a relation in $\eqref{eq:relation}$ with $k=3$:
\begin{equation}
 \alpha_{1} u_{1} + \alpha_{2} u_{2} + \alpha_{3} u_{3} = \beta_{1} v_{1} + \dots +  \beta_{m} v_{m}, 
 \label{eq:3-by-m-relation}
\end{equation}
where $\alpha_{1}$, $\alpha_{2}$, $\alpha_{3}$, $\beta_{1}$, \dots, $\beta_{m}$ are integers. 
Then one of the following two properties holds. 
\begin{enumerate}
 \item 
       $m=3$, $\alpha_1 = \alpha_2 = \alpha_3 = \beta_1 = \beta_2=\beta_3$ and 
       the relation \eqref{eq:3-by-m-relation} is a positive multiple of 
       \begin{eqnarray}
	\lefteqn{u_{\tri{a}{b_{1}}{b_{2}C}} + u_{\tri{a}{b_{2}}{b_{3}C}} + u_{\tri{a}{b_{3}}{b_{1}C}}}
	 & & \nonumber \\
	 & & \qquad = 
	 u_{\tri{a}{b_{2}}{b_{1}C}} + u_{\tri{a}{b_{3}}{b_{2}C}} + u_{\tri{a}{b_{1}}{b_{3}C}}
	 \label{eq:3-by-3-relation}
       \end{eqnarray}
 \item The relation \eqref{eq:3-by-m-relation} contains a relation of Theorem~\ref{thm:2:2}. 
\end{enumerate}
\end{theorem}

The statement in  (1) can not be strengthened because 
there exist $3$ by $m$ ($m\ge 3$) relations of type (2) such as
\begin{align*}
2 u_{\tri{a}{c}{\emptyset}}+  u_{\tri{a}{b}{c}} + u_{\tri{b}{c}{a}}
&=u_{\tri{a}{b}{\emptyset}} + u_{\tri{b}{c}{\emptyset}} + 2 u_{\tri{a}{c}{b}}\\
&(=u_{\tri{ab}{c}{\emptyset}} + u_{\tri{a}{bc}{\emptyset}}),
\end{align*}
and
\begin{align*}
&u_{\tri{a}{b_3}{b_1b_2}}+2 u_{\tri{a}{b_2}{b_1}} + u_{\tri{a}{b_1}{\emptyset}}\\
&\qquad =u_{\tri{a}{b_2}{b_1b_3}} + u_{\tri{a}{b_3}{b_1}} + u_{\tri{a}{b_1}{b_2}}+ u_{\tri{a}{b_2}{\emptyset}}.
\end{align*}

\begin{proof} 
Let $u=\alpha_{1}u_{1}+\alpha_{2}u_{2}+\alpha_{3}u_{3}$. 
From~(\ref{eq:relation-sum1}), we have 
\begin{equation}
 \alpha_1 + \alpha_2 + \alpha_3 = \beta_1 + \dots + \beta_m
  , \quad 
  \alpha_1 > 0, \alpha_2 > 0, \alpha_3 > 0, 
  \beta_1 > 0, \dots, \beta_m > 0. 
 \nonumber 
\end{equation}
Because $u_{1}, u_{2}, u_{3}, v_{1}, \dots, v_{m}$ are extreme rays of $\cK^*(N)$, 
we have $m \geq 3$. 
Furthermore, because $u_{1}, u_{2}, v_{1}, \dots, v_{m}$ are distinct elementary imsets, 
the following properties hold: 
\begin{equation}
 u_{1} \neq u_{2},\quad u_{1}, u_{2} \neq v_{1}, \dots, v_{m},\quad  v_{i} \neq v_{j}. 
 \label{eq:proof_theorem4:1}
\end{equation}
We write the elementary imsets in (\ref{eq:3-by-m-relation}) as 
$u_{1}=\utri{a_{1}}{b_{1}}{C_{1}}, u_{2}=\utri{a_{2}}{b_{2}}{C_{2}}, u_{3}=\utri{a_{3}}{b_{3}}{C_{3}}, v_{1}=\utri{a_{1}'}{b_{1}'}{C_{1}'}, \dots, v_{m}=\utri{a_{m}'}{b_{m}'}{C_{m}'}$. 
As in the proof of Theorem~\ref{thm:2:2}, we define the weight $w$ for a relation as follows:
\begin{equation}
 w = \sum_{ \{S \subseteq N \: \mid \: \abs{S}=\abs{C_{1}}+1 \} }(a_{1}u_{1}(S)+a_{2}u_{2}(S)+a_{3}u_{3}(S)). 
  \nonumber 
\end{equation}
Remember that the terms are ordered such that $u_{1} < u_{2} < u_{3}$ and $v_{1} < v_{2} < \dots < v_{m}$ 
according to the graded reverse lexicographic order of elementary imsets. 
Then the possible patterns are classified as follows:
\begin{enumerate}[(i)]
 \item $w=-2(\alpha_{1}+\alpha_{2}+\alpha_{3}) \Leftrightarrow \abs{C_{1}}=\abs{C_{2}}=\abs{C_{3}}$,
       \label{enum:thm4:1}
 \item $w=-2(\alpha_{1}+\alpha_{2})+\alpha_{3} \Leftrightarrow \abs{C_{1}}=\abs{C_{2}}=\abs{C_{3}}-1$,
       \label{enum:thm4:2}
 \item $w=-2\alpha_{1}+\alpha_{2}+\alpha_{3} \Leftrightarrow \abs{C_{1}}=\abs{C_{2}}-1=\abs{C_{3}}-1$,
       \label{enum:thm4:3}
 \item $w=-2\alpha_{1}+\alpha_{2} \Leftrightarrow \abs{C_{1}}=\abs{C_{2}}-1 \leq \abs{C_{3}}-2$.
       \label{enum:thm4:4}
 \item $w=-2\alpha_{1} \Leftrightarrow \abs{C_{1}} \leq \abs{C_{2}}-2 \leq \abs{C_{3}}-2$.  
\label{enum:thm4:5}
 \item $w=-2(\alpha_{1}+\alpha_{2}) \Leftrightarrow \abs{C_{1}} = \abs{C_{2}} \leq \abs{C_{3}}-2$ 
\label{enum:thm4:6}
\end{enumerate}

First, we consider the case \eqref{enum:thm4:5} of $w=-2\alpha_{1}$. 
Because $u_{1} \neq v_{i}$ from (52), 
there exists $i' \in \{1,\dots m\}$ 
such that $v_{i'}=\tri{a_{1}}{d'}{C_{1}}$, $d' \neq b_{1}$. 
To cancel out the value $v_{i'}(d'C_{1})=-\alpha_{i'} < 0$ at $d'C_{1}$, 
there exists $i'' \in \{1,\dots m\}$ 
such that $v_{i''}=\tri{a''}{b''}{d'C_{1}}$, $d' \neq b_{1}$. 
Furthermore there exists $d'' \in \{a'', b''\}$ such that $u_{i}(d'd''C_{1})=0$ for $i=1,2,3$. 
To cancel out the value $v_{i''}(d'd''C_{1})=-\alpha_{i''} < 0$ at $d'C_{1}$, 
there exists $v_{i'''}=\tri{a'''}{b'''}{d'd''C_{1}}$. 
In the same way, there exists $d'''\in\{a''', b'''\}$ such that $u_{i}(d'd''d'''C_{1})=0$ for $i=1,2,3$. 
Since it is impossible to cancel out the value $v_{i'''}(d'd''d'''C_{1}) = -\alpha_{i'''} < 0$ 
at $d'd''d'''C_{1}$, this case can be ignored. 

The case \eqref{enum:thm4:6} is almost the same as in the case of \eqref{enum:thm4:5}. 

We consider the case~\eqref{enum:thm4:4}. 
If $\abs{C_{2}}-1 < \abs{C_{3}}-2$, then it contradicts (\ref{eq:proof_theorem4:1}). 
Therefore, we consider the case of $\abs{C_{2}}-1 = \abs{C_{3}}-2$. 
In this case, there exist $m_{1}$ and $m_{2}$ such that 
\begin{eqnarray}
 1 \leq m_{1} \leq m_{2} \leq m, \quad  
 \alpha_1 = \beta_1 + \dots + \beta_{m_{1}} , \quad 
  C_{1} = C_{1}' = \dots = C_{m_{1}}', 
  \nonumber \\ 
 \alpha_3 = \beta_{m_{2}+1} + \dots + \beta_m , \quad
  \abC{3} = \abCd{m_{2}+1} = \dots = \abCd{m}, 
  \nonumber 
\end{eqnarray}
and we have $\abs{C_{2}}=\abs{C_{m_{1}+1}'}=\dots=\abs{C_{m_{2}}'}$. 
Now  let $\abC{1} = \{d_{1}, d_{2}\} \cup C_{1}$. 
From (\ref{eq:proof_theorem4:1}) and the same argument of the case \eqref{enum:thm4:5}, we can write $C_{2}=d_{1} C_{1}$ without loss of generality. 
Furthermore from (\ref{eq:proof_theorem4:1}), we have $\alpha_{1} \leq \alpha_{2}$. 
If $u(\abC{1}) \leq 0$, then $\alpha_{1}u_{1}+\alpha_{2}u_{2}$ has a relation of Theorem~\ref{thm:2:2}. 
Therefore, we have $u(\abC{1}) > 0$. 
In the same way, we have $u(C_{3}) > 0$. 
Hence there exist $i' \in \{m_{2}+1,\dots,m\}$ and $j' \in \{1,\dots,m_{1}\}$ such that 
$C_{i'}'=\abC{1}$ and $\abCd{j} \in C_{3}$. 
This means $\abC{1} \neq C_{3}$. 
Since $C_{1} \subseteq d_{1}C_{1} = C_{2} \subseteq \abC{1} = C_{i'}' \subseteq \abCd{i'} = \abC{3}$, 
we can write $\abC{3}=d_{2}d_{3}d_{4}C_{2}=d_{1}d_{2}d_{3}d_{4}C_{1}$. 
From (\ref{eq:proof_theorem4:1}), we have $d_{2}C_{1} \subseteq C_{3}$. 
This means that $C_{3}$ is equal to $d_{2}d_{3}C_{1}$ or $d_{2}d_{4}C_{1}$. 
Without loss of generality, we assume that $C_{3}=d_{2}d_{3}C_{1}$. 
In the same way, we can show that $\abC{2}$ is equal to $\abC{3}\setminus d_{1}$ or 
$\abC{3}\setminus d_{4}$.
Since $\abC{2} \supset C_{2} = d_{1}C_{1}$, we have $\abC{2}=\abC{3} \setminus d_{4}=d_{1}d_{2}d_{3}C_{1}$ 
and $u_{2}(d_{1}d_{2}C_{1})=-1$. 
However, from $\abC{1}=d_{1}d_{2}C_{1}$ and $a_{1} \leq a_{2}$, $u(\abC{1})=\alpha_{1}-\alpha_{2} \leq 0$. 
This is a contradiction. 

Next, we consider the case~\eqref{enum:thm4:3}. 
If $u(C_{2}), u(C_{3}) > 0$, then a contradiction follows from the proof of Theorem~\ref{thm:2:2}. 
Therefore, either $u(C_{2}) \leq 0$ or $u(C_{3}) \leq 0$ holds. 
Without loss of generality, we assume that $u(C_{2}) \leq 0$. 
From (\ref{eq:proof_theorem4:1}), we have $u(\abC{1}) \leq 0$. 
Therefore, we obtain
\begin{eqnarray}
 0 \geq u(\abC{1}) = \alpha_{1} + \alpha_{2}u_{2}(\abC{1}) + \alpha_{3}u_{3}(\abC{1}) 
  \nonumber \\ 
   \Longleftrightarrow
  \frac{\alpha_{2}}{\alpha_{1}}u_{2}(\abC{1}) + \frac{\alpha_{3}}{\alpha_{1}}u_{3}(\abC{1}) \leq -1. 
  \nonumber 
\end{eqnarray}
From this we have $u_{2}(\abC{1}) = -1$ or $u_{3}(\abC{1}) = -1$. 
However, if we assume $u_{2}(\abC{1}) = -1$, then $\alpha_{1}u_{1}+\alpha_{2}u_{2}$ has a relation of Theorem~\ref{thm:2:2}. 
Therefore, we have $u_{2}(\abC{1}) = 0$ and $u_{3}(\abC{1}) = -1$. 
From this, we have $C_{2} \subseteq \abC{1} \subseteq \abC{3}$. 
In this case, if $u(C_{3}) > 0$, then a contradiction follows from the proof of Theorem~\ref{thm:2:2}. 
Therefore $u(C_{3}) \leq 0$ holds. 
Furthermore, we have $u_{1}(C_{3}) = -1$. 
This means that $\alpha_{1}u_{1}+\alpha_{3}u_{3}$ has a relation of Theorem~\ref{thm:2:2}. 

The case~\eqref{enum:thm4:2} is almost the same as in the case of \eqref{enum:thm4:3}. 

Finally, we consider the case~\eqref{enum:thm4:1}. 
Since $C_{i}' \in \{C_{1},C_{2},C_{3}\}$, $\abCd{i} \in \{\abC{1},\abC{2},\abC{3}\}$ 
and (\ref{eq:proof_theorem4:1}), 
we have $m=3$ and $C_{i}'=C_{i}$. 
In this case, $\abCd{1} = \abC{2}$ or $\abCd{1} = \abC{3}$ holds. 
Without loss of generality, we assume $\abCd{1} = \abC{2}$. 
Then, we have $\abCd{2} = \abC{3}$, $\abCd{3} = \abC{1}$ 
and $\alpha_{1}=\alpha_{2}=\alpha_{3}=\beta_{1}=\beta_{2}=\beta_{3}$. 
Therefore the relation is a positive multiple of the three by three relation in (\ref{eq:3-by-3-relation}).
\end{proof}

\section{Some computational results on Markov basis for imsets}
\label{sec:markov-basis}


In this section we give some computational results on Markov basis for the configuration  $\cfg_{N}$.
Before giving them we discuss motivation for studying Markov basis 
from the viewpoint of inferences on conditional independence statements.
In Theorem \ref{thm:kernel-lattice-basis} we showed that the set of
moves in \eqref{eq:2-2-relation-kernel} contains a lattice basis for
$\ker_\ZZ \cfg_N$.  This means that if we allow negative coefficients,
then all relations with integer coefficients can be derived from the
two by two relations of Theorem \ref{thm:2:2}.
On the other hand, recall that the three by three relation in Theorem 
\ref{thm:3:3} can not be derived from two by two relations in the sense of our discussion after
\eqref{eq:3-3-equivalence}.  This is because in ``applying a rule'', we do not
allow subtracting a non-existing term from either side of a relation.
If a  Markov basis is available, then we can derive any relation without ever subtracting
a non-existing term starting from a given initial relation.  For definitions and details of Markov
basis refer to Chapter 5 of Sturmfels\cite{sturmfels1996} or
Chapter 1 of Drton, Sturmfels and Sullivant\cite{drton-sturmfels-sullivant-lecture}.

Let $\cG(\cfg_{N})$ be a (minimal) Markov basis for $\cfg_{N}$ 
and $\abs{\cG(\cfg_{N})}$ be the number of elements in $\cG(\cfg_{N})$. 
If $|N| \leq 2$, then there is no relation in $\cfg_{N}$ and hence $\ker\cfg_{N}=\emptyset$. 
If $|N|=3$, then there is only the relation of \eqref{eq:2:2-relation} with $C=\emptyset$. 
We now give some computational results on Markov basis for $N =4,5$. 
In the following results, 
we identify two relations if
one is obtained by permuting characters in another. 
For example, we identify the following two relations: 
$u_{\tri{a}{b}{d}} + u_{\tri{a}{c}{bd}}=u_{\tri{a}{c}{d}} + u_{\tri{a}{b}{cd}}$  
and 
$u_{\tri{a}{b}{c}} + u_{\tri{b}{d}{ac}}=u_{\tri{b}{d}{c}} + u_{\tri{a}{b}{cd}}$. 
Let $\cG^{R}(\cfg_{N})$ be the set of representative elements in $\cG(\cfg_{N})$ under this equivalence relation. 
Let $\abs{\cG^{R}(\cfg_{N})}$ be the number of elements in $\cG^{R}(\cfg_{N})$. 

Let $\cfg_{\tri{A}{B}{C}}$ be a sub-configuration of 
$\cfg_N$, $N=ABC$, 
which consists of all of the elementary imsets $u_{\tri{a}{b}{\Gamma}}$ 
such that $a \in A$, $b \in B$, $C \subseteq \Gamma \subseteq ABC$, $a < b$ and $\Gamma\cap a=\Gamma\cap b=\emptyset$. 
Markov basis for $\cfg_{\tri{A}{B}{C}}$ represents the relations which appear in 
the face $F_{\tri{A}{B}{C}}$ of $\cK^*(N)$.
Markov bases for the subconfiguration of elementary imsets belonging to
$F_{\tri{A}{B}{C}}$ is of some independent interest because they
generate a combinatorial  pure subring in the sense of Ohsugi, Herzog and {Hibi\cite{ 
ohsugi-herzog-hibi2000,ohsugi-2007,ohsugi-hibi2010}}.
Let $\cG^{R}({\tri{A}{B}{C}})$ be the representatives of Markov basis for $\cfg_{\tri{A}{B}{C}}$
and $\abs{\cG^{R}({\tri{A}{B}{C}})}$ be the number of elements in $\cG^{R}({\tri{A}{B}{C}})$. 
Let 
$\cG^{R}(\cT(N))
= \sum_{ABC=N} \cG^{R}({\tri{A}{B}{C}})$
and 
$\abs{\cG^{R}(\cT(N))}$ be the number of elements in 
${\cG^{R}(\cT(N))}$. 

%

For a relation of \eqref{eq:relation}, we call 
$\alpha_{1}+\dots+\alpha_{k}$ 
the degree. 
Note that $\sum_{i=1}^{k}\alpha_{i}=\sum_{j=1}^{m}\beta_{m}$. 
We give $\cG^{R}(\cfg_{N})$ and 
$\cG^{R}(\cT(N))$
at each degree for $|N|=4,5$ in 
Table~\ref{tab:computaional_results_n4} and \ref{tab:computaional_results_n5}. 
The Markov basis of $\cfg_{N}$ were computed with {4ti2\cite{4ti2}}. 
Since the computation has not completed 
when we used the graded reverse lexicographic order in (\ref{eq:grlex-for-elementary}), 
we give the result for the following order: 
\begin{eqnarray*}
u_{\tri{a}{b}{C}} \precw u_{\tri{a'}{b'}{C'}} \quad
\text{if}\quad
\begin{cases}
2^{C}>2^{C'}, & \\
\text{or }  C=C' \text{ and } b>b',&\\
\text{or }  C=C' \text{ and } b=b' \text{ and } a>a' , &
\end{cases}
\end{eqnarray*}
where $2^{C} = \sum_{i \in C}2^{i}$. 

\begin{table}[thbp]
 \tbl{Number of representatives of elements  of Markov basis for $|N|=4$}{\hspace*{10cm}}
 \label{tab:computaional_results_n4}
 \begin{center}
  \begin{tabular}[c]{ccc}
   degree
   & $\cG^{R}(\cfg_{N})$ 
   & $\cG^{R}(\cT(N))$ 
   \\ \hline 
   2 & 2 & 2 \\
    3 & 1 & 1 \\ 
   4 & 4 & 4 \\ 
  \end{tabular}
 \end{center}
\end{table}
\begin{table}[thbp]
 \tbl{Number of representatives of elements of Markov basis for $|N|=5$}{\hspace*{10cm}}
 \label{tab:computaional_results_n5}
 \begin{center}
  \begin{tabular}[c]{ccc}
   degree
   & $\cG^{R}(\cfg_{N})$ 
   & $\cG^{R}(\cT(N))$ 
   \\ \hline 
   2 &   3 & 3 \\
   3 &   2 & 2 \\
   4 &  11 & 11 \\
   5 &  18 & 16 \\
   6 & 210 & 162 \\
   7 & 384 & 36 \\
   8 & 364 & 38 \\
   9 &  90 & 0 \\
   10 & 220 & 0 \\
   11 &  16 & 0 \\
   12 &  63 & 0 \\
  \end{tabular}
 \end{center}
\end{table}

In Table~\ref{tab:computaional_results_n4} for $|N|=4$, 
all of the representatives in $\cG^{R}(\cfg_{N})$ appear in $\cG^{R}(\cT(N))$. 
However, in Table~\ref{tab:computaional_results_n5} for $|N|=5$, 
there are many representatives which appear in $\cG^{R}(\cfg_{N})$ but do not appear in $\cG^{R}(\cT(N))$. 
Especially there is no relation with degree more than 8 in $\cG^{R}(\cT(N))$. 
This means that the facial structure of $\cK^*(N)$ is quite complicated. 
Furthermore, through this computational study, 
though there are a lot of coefficients more than $1$ in the Markov basis for $\cG^{R}(\cT(N))$ for $|N|=5$, 
we confirmed that the Markov basis for $\cG^{R}({\tri{A}{B}{C}}),  \abs{A}+\abs{B}+\abs{C} \leq 5$
is square free i.e., the coefficients in the relations are all $1$. 
It is future work to confirm whether this is true or not for $|N| \geq 6$. 


\appendix

Here we present a generalization of Proposition \ref{prop:product-skeletal} to
convex cones on a direct product of two finite sets.

\begin{proposition}
\label{prop:generalized-GH}
Let $X,Y$ be finite sets and let 
$G\subseteq \{g:X\to\mathbb{R}^+\}$,
$H\subseteq \{h:Y\to\mathbb{R}^+\}$
be pointed convex cones in ${\mathbb R}^{X}$ and
${\mathbb R}^{Y}$ respectively.
Define  
\[
 M(G,H):=\{f:X\times Y\to\mathbb{R}^+\mid \forall x,\ f(x,\cdot)\in H,\ \forall y,\ f(\cdot,y)\in G\}.
\]
Then $M(G,H)$ is a convex cone.
If $g$ and $h$ are extreme rays of $G$ and $H$, then 
$g(x)h(y)$ is an extreme ray of $M(G,H)$.
\end{proposition}

\begin{proof}
 Let $g$ and $h$ be  extreme rays of  $G$ and $H$ and write 
 $g(x)h(y)=f_1(x,y)+f_2(x,y)$ ($f_1,f_2\in M(G,H)$).
 Note that there exist $x'$, $y'$ such that $g(x')>0$ and $h(y')>0$.
 For $x'$ with  $g(x')>0$ we have
 \begin{align*}
  h(y) = \frac{f_1(x',y)}{g(x')} + \frac{f_2(x',y)}{g(x')},\ \ \forall y.
 \end{align*}
 Since $h$ is an extreme ray of $H$ and 
 the two terms on the right-hand side belong to $H$ as functions of $y$ by the definition of 
  $M(G,H)$, $f_1(x',y)=c_1(x')h(y)$,  $f_2(x',y)=c_2(x')h(y)$ for some 
 $c_1(x') ,c_2(x')$.
 For $x$ with  $g(x)=0$, we have
 $f_1(x,y)=f_2(x,y)=0$ ($\forall y$) and 
 we can put  $c_1(x)=c_2(x)=0$. Then for all $x,y$
 we can write $f_1(x,y)=c_1(x)h(y),f_2(x,y)=c_2(x)h(y)$.
 In particular for $y'$ with  $h(y')>0$, 
 $c_1(x)=f_1(x,y')/h(y')$, $c_2(x)=f_2(x,y')/h(y')$
 belong to $G$ as functions of $x$ and 
  \begin{align*}
  g(x) = \frac{f_1(x,y')+f_2(x,y')}{h(y')} = c_1(x) + c_2(x).
 \end{align*}
 Hence  $c_1,c_2$ are proportional to $g$.
\end{proof}

We now show that Proposition \ref{prop:product-skeletal}  is a special case of 
Proposition \ref{prop:generalized-GH}.
We note that skeletal functions are also extreme rays of 
$\mathcal{K}_+(N):=\mathcal{K}(N)\cap \{f\mid f\geq 0\}$. 
In fact, if  $f\in \mathcal{K}_{\ell}(N)$ and 
$f=f_1+f_2$ ($f_1,f_2\in \mathcal{K}_+(N)$), then it is easily seen that
$f_1,f_2\in \mathcal{K}_{\ell}(N)$.


 Let $X=2^A,Y=2^B$, $G=\mathcal{K}_+(A)$,  $H=\mathcal{K}_+(B)$, 
 and let $g$ and $h$ be skeletal on $A$ and $B$.
 Then $g$ and $h$ are extreme rays of $G$ and $H$. 
 Hence by Proposition \ref{prop:generalized-GH} $f$ is an extreme ray of $M(G,H)$. 
 Now from the definition of supermodularity $\mathcal{K}_{\ell}(AB)\subseteq \mathcal{K}_+(AB)\subseteq M(G,H)$.  Therefore 
 $f$ is an extreme ray of $\mathcal{K}_{\ell}(AB)$.

\noindent
\begin{remark}
We can not put $G=\mathcal{K}_{\ell}(A),H=\mathcal{K}_{\ell}(B)$
in the above proof, because $\mathcal{K}_{\ell}(AB)\not\subseteq M(G,H)$.
For example, let $f(S)=1_{\{ab\subseteq S\}}$ ($a\in A,b\in B$), then
$f\in \mathcal{K}_{\ell}(AB)$ but 
$f(AS)=1_{\{b\in S\}}$ ($\forall S\subseteq B$) does not belong to 
$\mathcal{K}_{\ell}(B)$ and $f\notin M(\mathcal{K}_{\ell}(A),\mathcal{K}_{\ell}(B))$.
\end{remark}

\bigskip
\noindent
{\bf Acknowledgment.} \quad
We are grateful to Milan Studen\'y for very detailed and useful comments.

\bibliographystyle{ws-procs9x6}
\bibliography{imset}

\begin{thebibliography}{10}

\bibitem{stu2005}
M.~Studen\'y, {\em Probabilistic Conditional Independence Structures}
  (Springer-Verlag, London, 2005).

\bibitem{studeny-vomlel-hemmecke-2010}
M.~Studen{\'y}, J.~Vomlel and R.~Hemmecke, {\em Internat. J. Approx. Reason.}
  {\bf 51}, 573 (2010).

\bibitem{studeny-vomlel-2011}
M.~Studen{\'y} and J.~Vomlel, {\em Internat. J. Approx. Reason.} {\bf 52}, 627
  (2011).

\bibitem{grunbaum-2003}
B.~Gr{\"u}nbaum, {\em Convex Polytopes}, Graduate Texts in Mathematics,
  Vol.~221, second edn. (Springer-Verlag, New York, 2003).
\newblock Prepared and with a preface by Volker Kaibel, Victor Klee and
  G{\"u}nter M. Ziegler.

\bibitem{kashimura-takemura-2011}
T.~Kashimura and A.~Takemura, Standard imsets for undirected and chain
  graphical models arXiv:1102.2927v1, (2011).

\bibitem{sturmfels1996}
B.~Sturmfels, {\em Gr\"obner Bases and Convex Polytopes}, University Lecture
  Series, Vol.~8 (American Mathematical Society, Providence, RI, 1996).

\bibitem{kashimura-etal-2by2}
T.~Kashimura, T.~Sei, A.~Takemura and K.~Tanaka, {\em Journal of Algebraic
  Statistics} {\bf 2}, 14 (2011).

\bibitem{hemmecke-etal-2008}
R.~Hemmecke, J.~Morton, A.~Shiu, B.~Sturmfels and O.~Wienand, {\em Combin.
  Probab. Comput.} {\bf 17}, 239 (2008).

\bibitem{stu1994}
M.~Studen{\'y}, {\em International Journal of General System} {\bf 23}, 123
  (1994/1995).

\bibitem{studeny-etal-five-variables}
M.~Studen\'y, R.~R. Bouckaert and T.~Ko{\v c}ka, Extreme supermodular set
  functions over five variables Research report n.\ 1977, Institute of
  Information Theory and Automation, Prague, (2000).

\bibitem{topkis1998}
D.~M. Topkis, {\em Supermodularity and Complementarity} (Princeton University
  Press, Princeton, N.J., 1998).

\bibitem{rosenmuller-weidner-1974}
J.~Rosenm{\"u}ller and H.~G. Weidner, {\em Discrete Math.} {\bf 10}, 343
  (1974).

\bibitem{drton-sturmfels-sullivant-lecture}
M.~Drton, B.~Sturmfels and S.~Sullivant, {\em Lectures on algebraic
  statistics}, Oberwolfach Seminars, Vol.~39 (Birkh\"auser Verlag, Basel,
  2009).

\bibitem{ohsugi-herzog-hibi2000}
H.~Ohsugi, J.~Herzog and T.~Hibi, {\em Osaka J. Math.} {\bf 37}, 745 (2000).

\bibitem{ohsugi-2007}
H.~Ohsugi, {\em Comment. Math. Univ. St. Pauli} {\bf 56}, 27 (2007).

\bibitem{ohsugi-hibi2010}
H.~Ohsugi and T.~Hibi, {\em Ann. Inst. Statist. Math.} {\bf 62}, 639 (2010).

\bibitem{4ti2}
4ti2 team, 4ti2---a software package for algebraic, geometric and combinatorial
  problems on linear spaces {A}vailable at www.4ti2.de.

\end{thebibliography}
\end{document}